%% file: Laurain-Riviere.tex
\documentclass[a4paper,10pt]{article}
\usepackage{amssymb}
\usepackage[babel=true,kerning=true]{microtype}

\usepackage{amsthm,leqno}

\usepackage{amsfonts,amsmath,amssymb,enumerate}
\usepackage{graphicx}
\usepackage{tikz}
\usetikzlibrary{decorations.markings}
\usetikzlibrary{plotmarks}
\usetikzlibrary{patterns}

\newcommand{\ds}{\displaystyle}

\newcommand{\be}{\begin{equation*}}
\newcommand{\ee}{\end{equation*}}
\newcommand{\beq}{\begin{equation}}
\newcommand{\eeq}{\end{equation}}
\newcommand{\begincal}{\begin{eqnarray*}}
\newcommand{\fincal}{\end{eqnarray*}}
\newcommand{\lbeq}{\begin{leqno}}
\newcommand{\leeq}{\end{leqno}}

\newtheorem{thm}{Theorem}[section]
\newtheorem{lemma}{Lemma}[section]
\newtheorem{cor}{Corollary}[section]

\newtheorem{defi}{Definition}[section]

\newtheorem{rem}{Remark}[section]

\newcommand{\eps}{{\varepsilon}}
\newcommand{\R}{{\mathbb R}}
\newcommand{\C}{{\mathbb C}}
\newcommand{\N}{{\mathbb N}}
\newcommand{\h}{{\mathbb H}}
\newcommand{\Z}{{\mathbb Z}}

\def\ti{\widetilde}
\def\lf{\left}
\def\rg{\right}

\def\al{\alpha}
\def\la{\lambda}
\def\eps{\varepsilon}
\def\ds{\displaystyle}

\def\Om{\Omega}
\def\om{\omega}
\def\p{\partial}

\title{Angular Energy Quantization for Linear Elliptic Systems with Antisymmetric  Potentials and 
Applications}
\author{Paul Laurain \& Tristan Rivi\`ere}
\begin{document}
\maketitle
\begin{abstract} 
In the present work we establish a quantization result for the angular part of the energy of solutions to elliptic linear systems of Schr\"odinger type with antisymmetric potentials
in two dimension. This quantization is a consequence of uniform Lorentz-Wente type estimates  in degenerating annuli. We derive from this angular quantization the full energy quantization for general critical points to functionals which are conformally invariant or also for pseudo-holomorphic curves on degenerating Riemann surfaces.\\
\end{abstract}

\medskip

{\bf A.M.S. Classification :} 35J47, 35J20, 35J60, 58E20, 53C42, 49Q05, 53C21, 32Q65

\tableofcontents

\section*{Introduction}

 Conformal invariance is a fundamental property for many problems in physics and geometry. In the last decades it has become an important feature of many questions of 
 non-linear analysis
 too.  Elliptic conformally invariant lagrangians for instance share similar analysis behaviors : their Euler Lagrange equations are critical with respect to the function space
 naturally given by the lagrangian and, as a consequence, solutions to these Euler Lagrange equations are subject to {\it concentration compactness} phenomena. Questions such as the regularity of solutions or energy losses for sequences of solutions cannot be solved by robust general arguments in PDE but require instead a careful study of the interplay between
 the highest order part of the PDE and it's non-linearity.\\
For example, in dimension $2$, let $(\Sigma,h)$ be a closed Riemann surface, it has been proved,  see theorem I.2 of \cite{Ri3}, that every critical point of a conformally invariant functional, $u:\Sigma\rightarrow \R^n$, solves a system of the form\footnote{In coordinates this system reads 
$$
\forall i=1\cdots n\quad\quad-\Delta u_i = \sum_{j=1}^n\Omega^j_i \cdot \nabla u_j \mbox{ on } \Sigma
$$ 
where the $\cdot$  operation is the scalar product between the gradient vector fields $\nabla u_j$ and the different entries of the vector valued antisymmetric matrix $\Omega$.} 
\beq
\label{00}
-\Delta u = \Omega \cdot \nabla u \hbox{ on } \Sigma,
\eeq
where $\Omega\in so(n)\otimes T\Sigma$ and $\Delta$ is the negative Laplace-Beltrami operator  $\frac{1}{\sqrt{\vert h\vert}} \partial_i (\sqrt{\vert h\vert}h^{ij} \partial_j)$. The fundamental fact here that has been observed in \cite{Ri3} and exploited
in this work to obtain the H\"older continuity of $W^{1,2}-$solutions to (\ref{00}) is the {\bf anti-symmetry} of $\Omega$.

\medskip

 The analysis developed in \cite{Ri3} permitted to generalize to general 2-dimensional conformally invariant Lagrangians the use of {\it integrability by compensation theory}
 as it has been introduced originally by H.Wente in the framework of constant mean curvature immersions  in ${\mathbb R}^3$ solving the following {\it CMC-system}
\beq
\label{heq}
\Delta u = 2 u_x \wedge u_y  \hbox{ on } \Sigma.
\eeq
Solutions to this {\it CMC system} are in fact critical points to the following conformally invariant lagrangian
\[
E(u)=\frac{1}{2}\int_{\Sigma} |d u|_h^2\ dvol_h+\int_{\Sigma}u^\ast \omega
\]
where $\omega$ is a 2-form in ${\mathbb R}^3$ satisfying $d\om= 4\,dx_1\wedge dx_2\wedge dx_3$.
The natural space to consider the equation (\ref{heq}) is clearly the Sobolev space $W^{1,2}$.  The CMC-system (\ref{heq})
is critical for $W^{1,2}$ in the following sense : the r.h.s. of (\ref{heq}) is a-priori only in $L^1$. Classical Calderon Zygmund theory
tells us that derivatives of functions in $\Delta^{-1} L^1$ are in the weak $L^2$ space locally which is ''almost'' the information we started from.
Hence in a sense both the quadratic non-linearity for the gradient in the r.h.s of the system and the operator in the l.h.s. are at the same level from regularity point of view and it requires a more careful analysis in order to decide which one is leading the general dynamic of this system.

\medskip

H. Wente discovered the special role played by the jacobian in the r.h.s. of (\ref{heq}), see \cite{He} and references therein, and was able to prove that if $u$ satisfies (\ref{heq}) then 
\beq
\label{010}
\Vert \nabla u\Vert_{2}\leq C \Vert \nabla u\Vert_2^2\quad ,
\eeq
where $C$ is independent on $\Sigma$ and equals\footnote{This later fact has been discovered later on by Y.Ge in \cite{Ge}, see also \cite{He}.} $\sqrt{3/16\pi}$ . This inequality implies that if $\sqrt{3/16\pi} \|\nabla u\|_2<1$ then the solution is constant. This is what we call the {\it Bootstrap Test} and it is the key observation for proving Morrey estimates and deduce the H\"older regularity of general solutions to (\ref{heq}) which bootstraps easily in order to establish that solutions to (\ref{heq}) are in fact $C^\infty$ .

\medskip

Another analysis issue for this equation is to understand the behavior of sequences $u_k$ of solutions to the CMC system (\ref{heq}). Inequality (\ref{010}) tells us again that if the energy does not concentrate at a 
point then the system will behave locally like a linear system of the form $\Delta u=0$ : the non-linearity $2 u_x \wedge u_y$ in the r.h.s is dominated by the linear highest order term
$\Delta u$ in the l.h.s. . As a consequence of this fact we deduce that  sequences of solutions to (\ref{heq}) with uniformly bounded energy strongly converge in $C^p$ norm
for any $p\in {\mathbb N}$, modulo extraction of a subsequence and possibly away from finitely many points\footnote{In our notations we can have some $a_\infty^i$ that coincide with another.} in $\Sigma$, $\{a^1_\infty\cdots a^l_\infty\}$ where the $W^{1,2}-$norm concentrates, towards a smooth limit that solves also
(\ref{heq})
\[
u_k\longrightarrow u_\infty\quad\quad\mbox{ strongly in }C^p_{loc}(\Sigma\setminus\{a^1_\infty\cdots a^l_\infty\})\quad\forall p\in {\mathbb N}.
\]
 The question remains to understand how the convergence at the concentration points $a^i_\infty$ fails to be strong, in other words we want to understand how and how much energy
has been dissipated at the points $a^i_\infty$. A careful analysis shows that the energy is lost by concentrating solution on $\R^2$ of the CMC system (\ref{heq}), the so called {\it bubbles}, that converge to the $a_\infty^i$ : there exists points in $\Sigma$ $a_k^{i}\rightarrow a_\infty^i$ and a familly of sequences of radii $\la_k^i$ converging to zero such that
\[
u_k(\la_k^i x+a_k^i)\longrightarrow \omega^i(x) \mbox{ strongly in }C^p_{loc}(\R^2\setminus\{\mbox{finitely many points}\})\quad\forall p\in {\mathbb N}.
\]
where $\omega^i$ denote the bubbles, solutions on $\R^2$ of the CMC system (\ref{heq}). Because of the nature of the convergence it is clear that the Dirichlet energy lost
in the convergences amount at least to the sum of the Dirichlet energies of the {\it bubbles} $\omega^i$ :
\beq
\label{011}
\liminf_{k\rightarrow +\infty} \int_{\Sigma}|du_k|^2_{h}\ dvol_h\ge\int_{\Sigma}|du_\infty|^2_{h}\ dvol_h+\sum_{i=1}^l\int_{\R^2}|\nabla \omega^i|^2\ dx_1\,dx_2
\eeq
The question remains to understand if the inequality in (\ref{011}) is strict or is in fact an equality. This question for general conformally
invariant problems is known as the {\it energy quantization question : } is the loss of energy only concentrated in the forming {\it bubbles} or is
there any additional dissipation in the intermediate regions between the {\it bubbles} and shrinking at the limiting concentration points $a_\infty^i$
 in the so called {\it neck region}. Since the work of Sacks and Uhlenbeck \cite{SaU} where it has been maybe first considered, in the particular
  framework of minimizing harmonic maps from a Riemann surface into a manifold, this question has generated a special interest, intensive researches and
  several detailed results have been obtained in the last decades  on the subject. We refer to  \cite{Ri2} and reference therein for a survey on the 
  {\it energy quantization results}.
  Positive results establishing {\it energy quantization} (i.e. the inequality in (\ref{011}) is in fact an equality) often make use of some special \underbar{geometric} objects such as isoperimetric inequality or the hopf differential , see for instance \cite{Jo} or \cite{Pa}.
  In \cite{LR1} and \cite{LR2} the second author in collaboration with F.H. Lin introduced a more \underbar{functional analysis} type technique based on the use of the interpolation Lorentz spaces in order to prove {\it energy quantization results} in the special cases where
  the non-linearity of the conformally invariant PDE can be written as a linear combination of jacobians of $W^{1,2}-$functions. Using this technique we can for instance prove that equality holds in (\ref{011}) : {\it energy quantization} holds for the CMC-system,  the whole
  loss of energy exclusively arises in the bubbles. The main step in the proof consists in using an improvement of Wente inequality
  (\ref{010}) which has been obtained by L.Tartar and R.Coifman, P.L.Lions, Y.Meyer and S.Semmes in \cite{CLMS}. This improved {\it Lorentz-Wente type inequality} reads 
\beq
\label{012}
\Vert \nabla u\Vert_{L^{2,1}}\leq C\ \Vert \nabla u\Vert_2^2\quad ,
\eeq
 where this time $C$ depends a-priori on $(\Sigma,h)$ and where $L^{2,1}$ denotes the Lorentz space ''slightly'' smaller than $L^2$ given by the space of measurable
 function $f$ on $\Sigma$ satisfying
 \[
 \int_0^\infty \vert \{x\in\Omega\hbox{ ¬†s.t. } \vert f(x) \vert \geq \lambda \}\vert^\frac{1}{2}  \, d\lambda < +\infty\quad .
 \]
 The goal of the present paper is to extend energy quantization results to sequences of
 critical points to general conformally invariant lagrangians using functional analysis arguments in the style of \cite{LR2}.
 
 \medskip
 
 The constant in the inequality (\ref{012}) depends a-priori on the domain, at least on its conformal class since the equation is conformally invariant. But our {\it neck regions} 
 connecting the {\it bubbles} are conformally equivalent to degenerating annuli. The first task of the present work is to prove different lemma which give some uniform estimates on the $L^{2,1}$-norm of the gradient for solution to Wente type equations on degenerating annuli. This is the subject of  section \ref{wtl}. \\ 

In the following sections, we use these uniform estimates established in section~\ref{wtl}  for proving various quantization phenomena. In particular we get the quantization of the angular part of the gradient for solution of general elliptic second order systems with anti-symmetric potentials. What we mean here by the angular part is the component of the gradient in the orthogonal of the radial direction with respect to the nearest point of concentration. 
Precisely the first main result in the present work is the following.

%In fact we consider a modified distance in order to separate bubbles over bubbles. We also use some cut off functions to make our bubble well defined, see section \ref{qang} for precise definitions and notations.
\begin{thm}
\label{th-I.1}
Let $\Omega_k\in L^2(B_1,so(n) \otimes \R^2)$ and let $u_k\in W^{2,1}(B_1,\R^n)$ be a sequence of solutions of
\beq
\label{eq}
-\Delta u_k =\Omega_k\cdot \nabla u_k ,
\eeq 
with bounded energy, i.e.
\beq
\label{be}
\int_{B_1} \left(\vert \nabla u_k\vert^2 + \vert \Omega_k\vert^2 \right)\, dz  \leq M.
\eeq
Then there exists $\Omega_\infty \in L^2(B_1,so(n) \otimes \R^2)$ and $u_\infty\in W^{2, 1}(B_1,\R^n)$ a solution of $-\Delta u_\infty=\Omega_\infty\cdot \nabla u_\infty$ on $B_1$, 
 $l\in\N^*$ and
\begin{enumerate}
\item $\omega^1,\dots, \omega^l$ a family of solutions to system of the form
\[
-\Delta \om^i=\widetilde{\Om}^i_\infty\cdot\nabla\om^i\quad\quad\mbox{ on } \R^2
\]
where $\widetilde{\Om}^i_\infty\in L^2(\R^2,so(n)\otimes{\mathbb R}^2)$,
\item $a_k^1,\dots,a_k^l$ a family of converging sequences of points of $B_1$,
\item  $\lambda_k^1,\dots,\lambda_k^l$ a family of  sequences of positive reals converging all to zero,
\end{enumerate}
such that, up to a subsequence,  
$$\Omega_k \rightharpoonup \Omega_\infty \hbox{ in }  L^2_{loc}(B_1,so(n) \otimes \R^2),$$

$$u_k\rightarrow u_\infty \hbox{ on } W^{1,p}_{loc}(B_1 \setminus \{ a^1_\infty,\dots,a^l_\infty\}) \hbox{ for all } p\geq1  $$
and 
$$ \left\Vert \left\langle\nabla \left( u_k -u_\infty -\sum_{i=1}^l \omega^i_k \right),X_k\right\rangle \right\Vert_{L^2_{loc}(B_1)}\rightarrow 0,$$
where $\omega^i_k=\omega^i(a_k^i+\lambda_k^i \, . \, )$ and  $X_k= \nabla^\bot d_k$ with $\ds d_k= \min_{1\leq i\leq l } (\lambda_k^i + d(a_k^i,\, .\,))$.\\

Moreover, if we have $\Vert \Omega_k\Vert_\infty =O(1)$ or even just $\Omega_k =O(\vert \nabla u_k\vert)$ hence the convergence to the limit solution $u_\infty$ is in fact in $C^{1,\eta}_{loc}$ for all $\eta\in[0,1[$.
\end{thm}

The proof of theorem~\ref{th-I.1} is established through the iteration of the following result. It says that, if the $L^2$ norm of the potential $\Omega$ is below some threshold on every dyadic
sub-annuli of a given annulus,  the angular part of the Dirichlet energy of $u$ on a slightly smaller annulus is controlled by the maximal contribution of the Dirichlet energy of $u$
on the dyadic sub-annuli. Precisely we prove the following.
\begin{thm}
\label{lp}
There exists  $\delta>0$ such that for all  $r,R\in\R^*_+$ satisfying $2r<R$ for all $\Omega\in L^2(B_R\setminus B_r, so(n)\otimes\R^n)$ and $u\in W^{1,2}(B_R\setminus B_r, \R^n)$ satisfying
$$-\Delta u =\Omega \cdot \nabla u$$
and
$$ \sup_{r<\rho<\frac{R}{2}} \int_{B_{2\rho}\setminus B_\rho}\vert\Omega\vert^2 \, dz \leq \delta .$$
Then there exists $C>0$, independent of $u$, $r$ and $R$, such that 

$$\left\Vert \frac{1}{\rho} \frac{\p u}{\p \theta}\right\Vert^2_{L^{2}\left(B_{\frac{R}{2}}\setminus B_{2 r}\right)}\leq  C \Vert \nabla u\Vert_2 \left[\sup_{ r<\rho<\frac{R}{2}} \int_{B_{2\rho}\setminus B_\rho}\vert \nabla u\vert^2 \, dz \right]^{1/2}.$$
\end{thm}

Thanks to the quantization of the angular part for general elliptic systems with anti-symmetric potential, we can derive the {\it energy quantization} for critical points to an arbitrary  continuously conformally invariant elliptic Lagrangian with quadratic growth.
\begin{thm} 
\label{th-I.3}
 Let $N^k$  be a $C^2$ submanifold of $\R^m$ and $\omega$ be a $C^1$ $2$-form on $N^k$ such that the $L^\infty$-norm of $d\omega$ is bounded on $N^k$. Let $u_k$ be
 a sequence of critical points in $W^{1,2}(B_1,N^k)$ for the Lagrangian
\beq
\label{ci}
F(u) =\int_{B_1} \left[\vert \nabla u\vert^2 + \omega(u)(u_x,u_y)\right]\, dz
\eeq
with uniformly bounded energy, i.e.
\be
\Vert \nabla u_k\Vert_{2} \leq M.
\ee
Then there exists $ \Lambda \in C^0(TN\otimes \R^2, so(n)\otimes \R^2)$ and   $u_\infty\in W^{1, 2}(B_1,\R^n)$ a solution of $-\Delta u=\Lambda(u,\nabla u) \cdot \nabla u$ on $B_1$, 
 $l\in\N^*$ and
\begin{enumerate}
\item $\omega^1,\dots, \omega^l$ some non-constant $\Lambda$-bubbles, i.e non-constant solution of
$$-\Delta\omega = \Lambda(\omega, \nabla \omega)\cdot \nabla \omega \hbox{ on } \R^2,$$ 
\item $a^1_k,\dots,a^l_k$ a family of converging  sequences of points of $B_1$,
\item  $\lambda^1_k,\dots,\lambda^l_k$ a family of  sequences of positive reals converging all to zero,
\end{enumerate}
such that, up to a subsequence, 
$$u_k\rightarrow u_\infty \hbox{ on } C^{1,\eta}_{loc}(B_1\setminus \{ a^1_\infty,\dots,a^l_\infty\}) \hbox{ for all } \eta \in [0,1[$$
and 
$$ \left\Vert \nabla \left( u_k -u_\infty -\sum_{i=1}^l \omega^i_k \right)\right\Vert_{L^2_{loc}(B_1)}\rightarrow 0\quad,$$
where $\omega^i_k=\omega^i(a_k^i+\lambda_k^i \,.\,  )$.
\end{thm}
Previous works establishing {\it energy quantizations} for various conformally invariant elliptic Lagrangian usually require more regularity on the Lagrangian ( see for instance \cite{Jo},  \cite{Pa}, \cite{St}, \cite{DiT} , \cite{LiWa},\cite{Zhu}). For instance in \cite{Pa} or  \cite{LiWa} the {\it energy quantization} for harmonic maps
in two dimension is obtained through the application of the maximum principle to an ordinary differential inequality satisfied by the integration over concentric circles of the angular part  of the energy. The application of this procedure required an $L^\infty$ bound on the \underbar{derivatives of the second fundamental form}, see lemma 2.1 of \cite{LiWa}.
We insist on the fact that, in comparison to the previously existing energy quantization results, theorem~\ref{th-I.3} above requires an  $C^0$ bound on the \underbar{second fundamental form} only, which is a
weakening of the regularity assumption for the target of a magnitude one with respect to derivation. 
Another application of theorem~\ref{th-I.3} is the {\it energy quantization} for solutions to the prescribed mean curvature system, see corollary \ref{qh}, assuming only an $C^0$ bound on the mean curvature.  Again, previous {\it energy quantization} results were assuming uniform $C^1$ bounds on $H$, see \cite{BeRe} and \cite{CaMu}. Theorem~\ref{th-I.3} in the prescribed mean curvature system corresponds again for this problem to
weakening of the regularity assumption for the target of a magnitude one with respect to derivation in comparison to previous existing result .This weaker assumption moreover are the minimal one in order
for the Lagrangian to be continuously differentiable and this is why it coincides with the original one appearing in the formulation of the Heinz-Hildebrandt regularity conjecture in the 70's.
\\
 
In a last part, we present some more applications of the uniform Lorentz-Wente estimates we established in section 2. The first one, for instance, deals with sequences of 
pseudo holomorphic immersions of sequences of closed Riemann surfaces whose corresponding conformal class degenerate in the moduli space
of the underlying 2-dimensional manifold. In particular we give a new proof of the Gromov's compactness theorem in all generality, see theorem \ref{gr}. We also give some cohomological condition which garanties the {\it energy quantization} for sequences of harmonic maps on degenerating surfaces. Finally we give a very brief introduction to the quantization of the Willmore surface established recently in \cite{BR}, where the uniform  Lorentz-Wente estimates of section \ref{wtl} play a crucial role.

\medskip {\bf Acknowledgements}: The first author was visiting  the {\it Forschungsinstituts f\"ur Mathematik} at E.T.H. (Zurich)  when this work started, he would like to thank it for its hospitality and the excellent working conditions. The two authors would like to thank Francesca Da Lio for her useful comments on the manuscript. 

\section{Lorentz spaces and standard Wente's inequalities}
\label{lorentz}
Lorentz spaces seems to be the good spaces in order to get precise Wente's inequality, here we recall some classical facts about these spaces,  \cite{StWe} and \cite{Gra1} for details.\\

\begin{defi}  Let $D$ be a domain of $\R^k$, $p \in ]1,+\infty[$ and  $q \in [1,+\infty]$. The Lorentz space $L^{p,q}(D)$ is the set of measurable functions
$f : D\rightarrow \R$ such that
$$\Vert f\Vert_{p,q}= \left(\int_{0}^{+\infty} \left( t^\frac{1}{p} f^{**}(t)\right)^q \frac{dt}{t}\right)^\frac{1}{q} <+\infty \hbox{ if } q<+\infty$$
or
$$\Vert f\Vert_{p,\infty}= \sup \left(t^\frac{1}{p} f^{**}(t)\right ) \hbox{ if } q=+\infty$$
where $f^{**}(t)=\frac{1}{t}\int_0^tf^*(s)\, ds$ and   $f^*$ the decreasing rearrangement of $f$.
\end{defi}

Each $L^{p,q}$ may be seen as a deformation of $L^p$. For instance, we have
the strict inclusions
$$L^{p,1} \subset L^{p,q'} \subset L^{p,q''}\subset L^{p,\infty},$$
if $1 < q'< q''$. Moreover,
$$L^{p,p} = L^p.$$
Furthermore, if $\vert D \vert$ is finite, we have that for all $q$ and $q'$,
$$p > p' \Rightarrow  L^{p,q} \subset L^{p',q'}.$$
Finally, for $p \in ]1,+\infty[$ and  $q \in [1,+\infty]$, $L^{\frac{p}{p-1},\frac{q}{q-1}}$  is the dual of $L^{p,q}$.\\

In the case $p,q=2,1$ we can give an other definition. Let $\phi(\lambda)=\vert \{t\in [0,\vert D\vert] \hbox{ ¬†s.t. } f^*(t)\geq \lambda \}\vert$, we make the change of variable $t=\phi(\lambda)$ in the definition of the Lorentz-norm, which gives
$$\Vert f\Vert_{2,1}= 2 \int_{\sup\vert f \vert}^{0}  \phi^{-\frac{1}{2}}(\lambda) \lambda \phi'(\lambda) \,d\lambda .$$
Hence integrating by part, we get
\beq
\label{d1}
\Vert f\Vert_{2,1}= 4  \int^{+\infty}_{0}  \vert \{x\in\Omega\hbox{ ¬†s.t. } \vert f(x) \vert \geq \lambda \}\vert^\frac{1}{2}  \, d\lambda .
\eeq

To finish this preliminary, we quickly present the standard Wente's inequalities for elliptic system in Jacobian form. Indeed if $a$ and $b$ are in $W^{1,2}$ this is clear that $a_xb_y-a_yb_x$ is in $L^1$ but in fact thanks to its structure, it is subject to compensated phenomena and $a_xb_y-a_yb_x$ is in $\mathcal{H}^1$ the Hardy space which is a strict subspace of $L^1$ which as a better behaviour than $L^1$ with respect to Calderon-Zygmund theory, since the convolution of a function in $\mathcal{H}^1$  and the Green kernel $\log(\vert z\vert)$ is in $W^{2,1}$. This improvement of integrability is summarized in the following theorem.

\begin{lemma}[\cite{We},\cite{Tar},\cite{CLMS}]
\label{wente} 
Let  $a$ and $b$ be in $W^{1,2}(B_1)$. Let $\phi\in W^{1,1}_0(B_1)$ be the solution of
\be
\Delta \phi = a_xb_y-a_yb_x \hbox{ on } B_1\\
\ee
Then there exists a constant $C$ independent of  $\phi$ such that 
\beq
\label{w1}
\Vert \phi\Vert_\infty + \Vert \nabla \phi \Vert_{2,1} + \Vert \nabla^2 \phi \Vert_1 \leq C \Vert \nabla a\Vert \Vert \nabla b\Vert_2.
\eeq
\end{lemma}

A consequence of the previous theorem was obtain by Bethuel \cite{Bet} using a duality argument.
\begin{lemma} Let $a$ and $b$ be two measurable functions such that $\nabla a \in L^{2,\infty}(B_1)$ and $\nabla b \in L^{2}(B_1)$. Let $\phi\in W^{1,1}_0(B_1)$ be the solution of
\be
\Delta \phi = a_xb_y-a_yb_x \hbox{ on } B_1
\ee
Then there exists a constant $C$ independent of $\phi$ such that 
\beq
\label{w3}
\Vert \nabla \phi \Vert_{2} \leq C \Vert \nabla a\Vert_{2,\infty} \Vert \nabla b\Vert_2.
\eeq
\end{lemma}

{\bf Notation:} In the following, if we consider a norm with out specified its domains, it is implicitly assume that its domain of definition is the one of the function. We denote $B_R(p)$ the ball of radius $R$ centered at $p$ and we just denote $B_R$ when $p=0$.

\section{Wente type lemmas}
\label{wtl}

In this section we  are going to prove some uniform Wente's estimates on annuli whose conformal class is a priori not bounded. In fact those estimate were already known for the $L^\infty$-norm and the $L^2$-norm of the gradient, since it has been proved that the constant is in fact independent of the domain considered, see  \cite{To} and \cite{Ge}. But this fact is to our knowledge new for the $L^{2,1}$-norm of the gradient.

\begin{lemma} 
\label{l4}
let  $a,b \in W^{1,2}(B_1)$,  $0<\eps<\frac{1}{2}$,  and $\phi\in W^{1,1}_0(B_1\setminus B_\eps)$ a solution of 

$$ \Delta \phi = a_xb_y-a_yb_x \hbox{ on } B_1\setminus B_\eps .$$
Then $\nabla\phi\in L^{2,1}(B_1\setminus B_\eps)$ and, for all $\lambda >1$, there exists a positive constant $C(\lambda)$ independent of $\eps$ and $\phi$ such that 

$$\Vert \nabla \phi\Vert_{L^{2,1}(B_1\setminus B_{\lambda\eps} )} \leq C(\lambda) \Vert \nabla a\Vert_2\,  \Vert \nabla b\Vert_2 .$$
\end{lemma}

{\it Proof of lemma~\ref{l4}:}\\

First we consider a solution of our equation on the whole disk, that is to say $\varphi\in W^{1,1}_0(B_1)$ which satisfies
$$ \Delta \varphi = a_xb_y-a_yb_x \hbox{ on } B_1 .$$
Then thanks to the classical Wente's inequality (\ref{w1}), we have 
\beq
\label{21}
\Vert \varphi\Vert_{\infty} +  \Vert \nabla \varphi \Vert_{2,1} \leq C \Vert \nabla a\Vert_2\,  \Vert \nabla b\Vert_2 .,
\eeq
where $C$ is a positive constant independent of $\varphi$.\\

Then we set $\psi=\phi-\varphi$, which satisfies
\be
\begin{cases}
&\Delta \psi =0 \hbox{ on }  B_1\setminus  B_\eps,\\
& \psi=0  \hbox{ on }  \partial B_1,\\
&\psi=-\varphi  \hbox{ on }  \partial B_\eps.
\end{cases}
\ee
Hence $\widetilde{\psi}=\psi -\left(\int_{\partial B_\eps}  \psi \, d\sigma\right) \frac{\ln(\vert z\vert)}{2\pi\eps\ln(\eps)}$ satisfies the hypothesis of the lemma \ref{l1}, then 
$$ \Vert \nabla \widetilde{\psi}\Vert_{L^{2,1}(B_1\setminus B_{\lambda\eps})} \leq C(\lambda) \Vert \nabla \widetilde{\psi}\Vert_2  \hbox{ for all } \lambda>1. $$
Hence, computing the $L^2$-norm of the gradient of the logarithm on $B_1\setminus B_{\lambda\eps}$, we get that
\beq
\label{a} 
\Vert \nabla \widetilde{\psi}\Vert_{L^{2,1}(B_1\setminus B_{\lambda\eps})} \leq C(\lambda)\left(\Vert \nabla \psi \Vert_2 + \left(\int_{\partial B _\eps}  \vert \psi\vert \, d\sigma\right) \frac{1}{\eps\sqrt{\ln\left(\frac{1}{\eps}\right)}}\right).
\eeq
But $\psi$ is the harmonic on $B_1\setminus  B_\eps$ and is equal to $-\varphi$ on the boundary, then
\beq
\label{a0}
\Vert \nabla \psi\Vert_2\leq \Vert \nabla \varphi\Vert_2 \hbox{ and }\Vert \psi \Vert_\infty\leq \Vert \varphi \Vert_\infty.
\eeq
Hence we get  that 
\beq
\label{aa}
 \int_{\partial B_\eps}  \vert \psi\vert \, d\sigma\leq \eps  C(\lambda)   \Vert \nabla a\Vert_2\,  \Vert \nabla b\Vert_2 ,
 \eeq
which gives, using (\ref{a}) and (\ref{a0}), that
\beq
\label{a1}
\Vert \nabla \widetilde{\psi}\Vert_{L^{2,1}(B_1\setminus B_{\lambda\eps})} \leq  C(\lambda) \Vert \nabla a\Vert_2\,  \Vert \nabla b\Vert_2.
\eeq 
Finally, computing the $L^{2,1}$-norm of the gradient of the logarithm on $B_1\setminus B_{\lambda\eps}$, we get that
\beq
\label{a2}
\Vert \nabla  \ln \Vert_{L^{2,1}(B_1\setminus B_{\lambda\eps})} =4 \sqrt{\pi} \ln\left(\frac{1}{\lambda\eps}\right).
\eeq
Hence, thanks to (\ref{aa}), (\ref{a1}) and (\ref{a2}), we get that 
\beq
\label{a3}
\Vert \nabla \psi  \Vert_{L^{2,1}(B_1\setminus B_{\lambda\eps})}  \leq C(\lambda)\Vert \nabla a\Vert_2\,  \Vert \nabla b\Vert_2 . 
\eeq
Then, thanks to (\ref{21}) and (\ref{a3}), we get the desired estimate.\hfill$\square$

\begin{lemma}
\label{lm-2.2}
let  $a,b \in W^{1,2}(B_1)$,  $0<\eps<\frac{1}{4}$,  and $\phi\in W^{1,1}(B_1\setminus B_\eps)$ a solution of 
\beq
\label{lbis}
\left\{
\begin{split}
& \Delta \phi = a_xb_y-a_yb_x \hbox{ on } B_1\setminus B_\eps \\[5mm]
&\int_{\partial B_\eps} \phi\;d\sigma =0,\\[5mm]
&\left\vert \int_{\partial B_1 } \phi  \, d\sigma \right\vert \leq K,\\
\end{split}
\right.
\eeq
where $K$ is a constant independent of $\eps$. Then, for all $ 0<\lambda <1$,  there exists a positive constant $C(\lambda)$ independent of $\eps$ such that 

$$\Vert \nabla \phi\Vert_{L^{2,1}\left( B_{\lambda}\setminus   B_{\lambda^{-1}\eps }\right)} \leq C(\lambda)( \Vert \nabla a\Vert_2\,  \Vert \nabla b\Vert_2 + \Vert \nabla \phi \Vert_2 + 1)\quad .$$
\hfill $\Box$
\end{lemma}

\noindent{\it Proof of lemma~\ref{lm-2.2} :}\\

Let $u\in W^{1,1}(B_1\setminus B_\eps)$ be the solution of 
\be
\left\{
\begin{split}
& \Delta u= 0 \hbox{ on } B_1 \setminus B_\eps, \\[5mm]
&u=\phi \hbox{ on ¬†} \partial B_1 \cup \partial B_\eps.
\end{split}
\right.
\ee

Hence $\Vert \nabla u \Vert_2 \leq \Vert \nabla \phi \Vert_2$. Moreover thanks to lemma \ref{l3} and lemma \ref{l4} we have $\nabla u \in L^{2,1}( B_{\lambda} \setminus  B_{\lambda^{-1}\eps })$ and   $\nabla (u - \phi) \in L^{2,1}(B_1 \setminus  B_{\lambda^{-1}\eps })$ with

\be
\begin{split}
&\Vert \nabla u \Vert _{L^{2,1}( B_{\lambda} \setminus  B_{\lambda^{-1}\eps })} \leq C(\lambda)\left(\Vert \nabla \phi \Vert _{2} +1\right)\\[5mm]
&\Vert \nabla (u-\phi) \Vert _{L^{2,1}( B_{\lambda} \setminus  B_{\lambda^{-1}\eps })}  \leq C(\lambda) \Vert \nabla a \Vert _{2} \Vert \nabla b \Vert _{2},
\end{split}
\ee
 which proves lemma~\ref{lm-2.2}.\hfill$\square$\\

\noindent{\bf Remark:} As in lemma \ref{l3}  we cannot  control the $L^{2,1}$-norm of $\nabla \phi$ by its $L^{2}$-norm, as it is shown by the following example
$$z\mapsto \frac{ln\left(\frac{\vert z\vert}{\eps}\right)}{ln\left(\frac{1}{\eps}\right)}.$$

\begin{lemma}
\label{LR0}
let  $a,b \in W^{1,2}(B_1)$,  $0<\eps<\frac{1}{4}$ and  $\phi\in W^{1,2}(B_1\setminus B_\eps)$ be a solution of 
\beq
\Delta \phi = a_xb_y-a_yb_x \hbox{ on } B_1\setminus B_\eps .
\eeq
Moreover, we assume  that 
\beq
\label{y1}
\Vert \phi \Vert_{\infty}< +\infty.
\eeq
Then, for $0 < \lambda <1$ a positive constant $C(\lambda)$ independent of $\eps$ and $\phi$ such that 
\beq
\label{y2}
\ds\Vert \nabla \phi \Vert_{L^{2,1}\left(B_\lambda\setminus  B_{\lambda^{-1}\eps }\right)} \leq C(\lambda) \left( \ \Vert \nabla a\Vert_{{2}}\,  \Vert \nabla b\Vert_{2} + \|\phi\|_\infty \right).
\eeq
\end{lemma}

\noindent{\it Proof of lemma~\ref{LR0}. :}\\

We introduce first $\varphi\in W^{1,2}_0(B_1\setminus B_\eps)$ to be the unique solution to

$$
\left\{
\begin{array}{l}
\Delta \varphi = a_xb_y-a_yb_x \hbox{ on } B_1\setminus B_\eps \\[5mm]
\varphi=0 \hbox{ on }\p B_1\cup \p B_\eps.
\end{array}
\right.
$$
Then thanks to lemma \ref{l4}, we have 
\be
\Vert\nabla \varphi\Vert_{L^{2,1}(B_1\setminus B_{\lambda^{-1}\eps } )} \leq C(\lambda)\ \Vert \nabla a\Vert_{2}\,  \Vert \nabla b\Vert_2 ,
\ee
where $C(\lambda)$ is a positive constant depending on $\lambda$ but not on $\phi$ and $\eps$.\\

Then we set $\psi=\phi-\varphi$, which is harmonic. Thanks to standard estimates on harmonic function, see \cite{HaLi} for instance, there exists $C(\lambda)>0$ a positive constant independent of $\psi$ and $\eps$ such that 

$$ \Vert \psi \Vert_{L^{2,1}\left(B_\lambda \setminus  B_{\lambda^{-1}\eps }\right)} \leq C(\lambda) \Vert \psi\Vert_{L^{\infty}(\partial B_1 \cup \partial B_\eps)} \leq C  \Vert \phi\Vert_{L^{\infty}}.$$
Which proves the desired inequality and lemma~\ref{LR0} is proved.\hfill$\square$

\begin{lemma}
\label{LR1}
let  $a,b \in L^2(B_1)$,  $0<\eps<\frac{1}{4}$,  assume that $\nabla a\in L^{2,\infty}(B_1)$ and that $\nabla b\in L^2(B_1)$,  let $\phi\in W^{1,(2,\infty)}(B_1\setminus B_\eps)$ a solution of 
\beq
\Delta \phi = a_xb_y-a_yb_x \hbox{ on } B_1\setminus B_\eps ,\\
\eeq
Denote, for $\eps\le r\le1$, $\phi_0(r):=(2\pi\, r)^{-1}\,\int_{\p B_r(0)}\phi\, d\sigma$ and assume
\beq
\label{x1}
\int_{\eps}^1|\dot{\phi_0}|^2\ r\ dr<+\infty\quad.
\eeq
Then, for $0\leq \lambda<1$, there exists a positive constant $C(\lambda)>0$ independent of $\eps$ and $\phi$ such that 
\beq
\label{x2}
\begin{split}
\ds\Vert \nabla \phi \Vert_{L^{2}\left(B_\lambda\setminus  B_{\lambda^{-1}\eps }\right)} &\leq C(\lambda) \left( \ \Vert \nabla a\Vert_{{2,\infty}}\,  \Vert \nabla b\Vert_{2} + \|\nabla\phi_0\|_{L^2(B_1\setminus B_\eps)} \right. \\
& \left.+ \ \Vert\nabla \phi\Vert_{L^{2,\infty}(B_1\setminus B_\eps)}\right).
\end{split}
\eeq
\end{lemma}

\noindent {\it Proof of lemma~\ref{LR1} :}\\

First we consider  $\varphi\in W^{1,2}_0(B_1)$ to be the solution of

$$
 \lf\{
\begin{array}{l}
\Delta \varphi = a_xb_y-a_yb_x \hbox{ on } B_1 \\[5mm]
\varphi=0 \quad\hbox{ on }\p B_1
\end{array}
\rg.
.$$
Then thanks to the generalized Wente's inequality, see (\ref{w3}), we have 
\beq
\label{x4}
\Vert\nabla \varphi\Vert_2 \leq C\ \Vert \nabla a\Vert_{2,\infty}\,  \Vert \nabla b\Vert_2 \quad.
\eeq
Consider the difference $v:=\phi-\varphi-(\phi_0-\varphi_0)$, it is an harmonic function on $B_1\setminus B_\eps$ which does not have
$0-$frequency Fourier modes :

\[
v=\sum_{n\in{\Z}^\ast} (c_n \rho^n+d_n \rho^{-n})\ e^{i\,n\theta}
\]
which implies in particular that 
\beq
\label{x5}
\int_{\p B_\rho}\frac{\p v}{\p \nu}\, d\sigma=0 \hbox { for all }  \eps<\rho<1.
\eeq
Moreover, due to the assumption (\ref{x1}) and due to (\ref{x4}) we have

\beq
\label{x6}
\begin{split}
\|\nabla v\|_{L^{2,\infty}(B_1\setminus B_\eps)} &\le 2 \|\nabla\varphi\|_2+ \|\nabla\phi_0\|_2 +\Vert\nabla \phi\Vert_{L^{2,\infty}(B_1\setminus B_\eps)}  \\[5mm]
&\le C \left(\ \Vert \nabla a\Vert_{2,\infty}\,  \Vert \nabla b\Vert_2 +\|\nabla\phi_0\|_2 + \Vert\nabla \phi\Vert_{L^{2,\infty}(B_1\setminus B_\eps)}\right)\quad.
\end{split}
\eeq
Let $\lambda\in]0,1[$, then standard elliptic estimates on harmonic function give  that $\forall \rho\in (\lambda^{-1}\eps,\lambda)$
\beq
\label{x7}
\begin{array}{l}
\ds\|\nabla v\|_{L^\infty(\p B_\rho)}\le C(\lambda)\ \rho^{-1}\ \|\nabla v\|_{L^{2,\infty}(B_{\lambda^{-1}\rho}\setminus B_{\lambda\rho})}\\[5mm]
\ds\quad\quad\quad\le C(\lambda)\ \rho^{-1}\ \left(\Vert \nabla a\Vert_{2,\infty}\,  \Vert \nabla b\Vert_2 +\|\nabla\phi_0\|_2  + \Vert\nabla \phi\Vert_{L^{2,\infty}(B_1\setminus B_\eps)}\right)\quad.
\end{array}
\eeq
Denote $\Om_\eps:=B_\lambda\setminus B_{\lambda^{-1}\eps}$. We have that
\beq
\label{x7a}
\|\nabla v\|_{L^2(\Om_\eps)}=\sup_{\{X\ ; \ \|X\|_{L^2(\Om_\eps)}\le 1\}}\int_{\Om_\eps}\nabla v\cdot X \, dz
\eeq
 For such an $X\in L^2(\Om_\eps)$ we denote $\ti{X}$ it's extension by 0 in the complement of $\Om_\eps$ in $B_1$. Let $g$ be the solution 
 of
 \[
 \lf\{
 \begin{array}{l}
 \ds\Delta g=-div \ti{X}^\perp\quad\quad\mbox{ in }B_1\\[5mm]
 \ds g=0\quad\quad\mbox{ on  }\p B_1
 \end{array}
 \rg.
 \]
where $\ti{X}^\bot=(-\ti{X}_2,\ti{X}_1)$. We easily see that 
\beq
\label{x8}
\|\nabla g\|_{L^2(B_1)}\le C\ \|\ti{X}\|_{L^2(B_1)}\le C\quad.
\eeq
Poincar\'e lemma gives the existence of $f\in W^{1,2}(B_1)$ such that
\[
\ti{X}=\nabla f+\nabla^\perp g\quad.
\]
and we have
\beq
\label{x9}
\|\nabla f\|_{L^2(B_1)}\le \|\nabla g\|_{L^2(B_1)}+ \|\ti{X}\|_{L^2(B_1)}+\le C+1\quad.
\eeq
We have
\[
\int_{\Om_\eps}\nabla v\cdot X\, dz=\int_{\Om_\eps}\nabla v\cdot\nabla f\, dz+\int_{\Om_\eps}\nabla v\cdot\nabla^\perp g\, dz
\]
We write
\beq
\label{x10}
\begin{array}{l}
\ds\int_{\Om_\eps}\nabla v\cdot\nabla^\perp g\, dz=\int_{\p B_\lambda}\p_\tau v\ g\, d\sigma-\int_{\p B_{\lambda^{-1}\eps}}\p_\tau v\ g\, d\sigma \\[5mm]
\ds\quad\quad=\int_{\p B_\lambda}\p_\tau v\ (g-g_{\lambda})\, d\sigma-\int_{\p B_{\lambda^{-1}\eps}}\p_\tau v\ (g-g_{\lambda^{-1}\eps})\, d\sigma
\end{array}
\eeq
where $\p_\tau$ is the tangential derivative along the circles $\p  B_{\lambda}$ and $\p  B_{\lambda^{-1}\eps}$ and $g_{\lambda}$ (resp. $g_{\lambda^{-1}}$) denote the average of $g$ on $\p  B_{\lambda}$ (resp. $ B_{\lambda^{-1}\eps}$).

We have for any $\rho\in(0,1)$

\beq
\label{x11}
\frac{1}{\rho}\int_{\p B_\rho}|g-g_\rho| \, d\sigma \le C\ \|g\|_{H^{1/2}(\p B_{\rho})}\le C\ \|\nabla g\|_2\le \ C
\eeq

where $C$ is independent of $\rho$. Combining (\ref{x7}), (\ref{x11}) and (\ref{x10}) give in one hand

\beq
\label{x12}
\lf|\int_{\Om_\eps}\nabla v\cdot\nabla^\perp g\rg| \, dz \le C(\lambda)\ \|\nabla v\|_{L^{2,\infty}(B_1\setminus B_\eps)}\quad.
\eeq

In the other hand one using the fact that $v$ is harmonic and satisfies (\ref{x5}) we have

\beq
\label{x13}
\begin{array}{l}
\ds\int_{\Om_\eps}\nabla v\cdot\nabla f=\int_{\p  B_{\lambda}}\p_\nu v\ f-\int_{\p  B_{\lambda^{-1}\eps}}\p_\nu v\ f\\[5mm]
\ds\quad\quad=\int_{\p  B_{\lambda}}\p_\nu v\ (f-f_{\lambda})-\int_{\p  B_{\lambda^{-1}\eps}}\p_\nu v\ (f-f_{\lambda^{-1}\eps})
\end{array}
\eeq

We have for any $\rho\in(0,1)$

\beq
\label{x14}
\frac{1}{\rho}\int_{\p B_{\rho}}|f-f_\rho|\le C\ \|f\|_{H^{1/2}(\p B_{\rho})}\le C\ \|\nabla f\|_2\le \ C
\eeq

Combining now (\ref{x7}), (\ref{x13}) together with (\ref{x14}) we obtain

\beq
\label{x15}
\lf|\int_{\Om_\eps}\nabla v\cdot\nabla f\rg|\le C(\lambda) \ \|\nabla v\|_{L^{2,\infty}(B_1\setminus B_\eps)}\quad.
\eeq

Combining (\ref{x12}), (\ref{x15}) with (\ref{x7a}) gives

\beq
\label{x16}
\|\nabla v\|_{L^2(\Om_\eps)}\le C(\lambda)\ \|\nabla v\|_{L^{2,\infty}(B_1\setminus B_\eps)}
\eeq

This inequality together with (\ref{x1}), (\ref{x4}) gives (\ref{x2}) and the lemma is proved.\hfill$\square$

\section{Angular Energy Quantization for solutions to elliptic systems with anti-symmetric potential }
\label{qang}
The aim of this section is to prove that the angular part of the gradient of a bounded sequence of solutions of an elliptic system with anti-symmetric potential is always quantified. But before starting the proof of  the quantization, we remind some fact about elliptic systems with antisymmetric potential which have intensively studied by  the second author \cite{Ri3}.\\

Let  $\Omega\in L^2(B_1, so(n)\otimes\R^n)$ we consider  $u\in W^{1,2}(B_1, \R^n)$  a solution of the following  equation 
$$-\Delta u = \Omega\cdot \nabla u \hbox{ on } B_1 .$$

One of the fundamental fact  about this system is the discover a conservation law using a Coulomb gauge for $\Omega$ when its $L^2$-norm is small enough which is  the aim of the following theorem.
\begin{thm}[Theorem I.4 \cite{Ri3}]
\label{lx}
 There exists  $\eps_0>0$  such that for all $\Omega\in L^2(B_1, so(n)\otimes\R^2)$ satisfying
$$\int_{B_1} \vert \Omega\vert^2 \, dz\leq \eps_0, $$
then there exists $A\in W^{1,2}\cap L^\infty(B_1, Gl_n(\R))$ such that 
$$div(\nabla A -A\Omega)=0$$
and
$$ \int_{B_1}( \vert \nabla A\vert^2 +  \vert \nabla A^{-1}\vert^2 ) \, dz + dist(\{A,A^{-1}\},SO(n))  \leq C \int_{B_1} \vert \Omega \vert^2 \, dz ,$$
where $C$ is a constant independent of $\Omega$.
\end{thm}

Then, using this theorem and Poincar\'e's lemma, we get the existence  of $B\in W^{1,2}(B_1,M_n(\R))$ such that 
$$ div(A\,\nabla u)= \nabla^\bot B\cdot \nabla u $$
and 
$$ \int_{B_1} \vert \nabla B\vert^2 \, dz  \leq C \int_{B_1} \vert \Omega \vert^2 \, dz .$$
Hence the system is rewrite in Jacobian form and we can use standard Wente's estimates. In particular, this permits to prove three fundamental properties of the solutions of this equation which are the $\eps$-regularity, the energy gap for solution defined on the whole plane and the passage to the weak limit in the equation. This properties are summarized in the following theorem.

\begin{thm}
\label{ereg}  \cite{Ri3}, \cite{Ri6}
 There exists $\eps_0 >0$ and $C_q>0$, depending only on $q
 \in \N^*$,  such that if  $\Omega\in L^2(B_1, so(n)\otimes\R^2)$ (reps.  $L^2(\R^2, so(n)\otimes\R^2)$)  satisfies 
$$\Vert \Omega\Vert^2_2  \leq \eps_0, $$
 then
\begin{enumerate}

\item ($\eps$-regularity) If $u\in W^{1,2}(B_1,\R^n)$  satisfies
\beq 
-\Delta u = \Omega\cdot \nabla u \hbox{ on } B_1
\eeq
then we have the following estimate
$$ \Vert \nabla u\Vert_{L^q\left(B_{\frac{1}{4}}\right)}  \leq C_q \Vert \nabla u\Vert_2 \hbox{ for all } q\in \N^*.$$
\item (Energy gap)  If   $u\in W^{1,2}(\R^2,\R^n)$  satisfies
\beq 
-\Delta u = \Omega\cdot \nabla u \hbox{ on } \R^2
\eeq
then it is constant.
\item (Weak limit property) Let $\Omega_k\in L^2(B_1,so(n)\otimes \R^2)$ such that $\Omega_k$ weakly converge in $L^2$ to $\Omega$ and $u_k$ a bounded sequence in  $W^{1,2}(B_1,\R^n)$  which satisfies
$$-\Delta u_k = \Omega_k\cdot \nabla u_k \hbox{ on } B_1 .$$
Then, there exists a subsequence of $u_k$ which weakly converge in  $W^{1,2}(B_1,\R^n)$ to a solution of 
\beq 
-\Delta u = \Omega\cdot \nabla u \hbox{ on } B_1.
\eeq
\end{enumerate}
\end{thm}

\medskip

For the convenience of the reader we recall the arguments developed in \cite{Ri3} and \cite{Ri6} to prove theorem~\ref{ereg}.

\medskip 

\noindent{\it Proof of theorem~\ref{ereg} :}\\

In order to prove the $\eps$-regularity, let us prove that, for $\alpha >0$,  we have
 \beq
 \label{Mo}
\sup_{p\in B_{1/2}\;,\; 0<\rho<\frac{1}{2}}\rho^{-\al}\ \int_{B_\rho(p)}|\Delta u|\, dz \leq C\Vert \nabla u\Vert_{L^2(B_1)} \quad .
\eeq
A classical estimate on Riesz potentials gives
$$
|\nabla u|(p)\le C\frac{1}{|x|}\ast \chi_{B_{1/2}}\ |\Delta u|+C \Vert \nabla u\Vert_{L^2(B_1)}\qquad\forall\:\:p\in B_{1/4}\quad,
$$
where $\chi_{B_{1/2}}$ is the characteristic function of the ball $B_{\frac{1}{2}}$. Together with injections proved by Adams in \cite{Ad}, the latter shows that
$$
\Vert \nabla u \Vert_{L^{r}\left(B_{\frac{1}{4}}\right)} \leq C \Vert \nabla u\Vert_{L^2(B_1)}\quad ,
$$
for some $r>1$. Then bootstrapping this estimate, see lemma IV.1 of \cite{Ri6} or theorem 1.1 of \cite{ShTo}, we get 
$$\Vert  \nabla u\Vert_{L^q(B_{\frac{1}{4}})} \leq C_q \Vert \nabla u\Vert_{L^2(B_1)} \hbox{ for all } q\in \N^*,$$
which will prove the $\eps$-regularity.\\ 

In order to prove (\ref{Mo}), we assume that $\eps_0$ is small enough to apply theorem \ref{lx} . Hence there exists $A\in W^{1,2}\cap L^\infty(B_1, Gl_n(\R))$ and $B\in W^{1,2}\cap L^\infty(B_1, M_n(\R))$such that 

$$ \int_{B_1} ( \vert \nabla A\vert^2 + \vert \nabla B \vert^2) \, dz + dist(\{A,A^{-1}\},SO(n)) \leq C \int_{B_1} \vert \Omega \vert^2 \, dz .$$
and 
\beq
\begin{cases}
& div (A\,\nabla u) =\nabla^\bot B\cdot \nabla u ,\\[5mm]
& curl (A\,\nabla u) =\nabla^\bot A\cdot \nabla u .
 \end{cases}
 \eeq
 Let $p\in B_{\frac{1}{2}}$ and $0<\rho<\frac{1}{2}$, we proceed by introducing on $B_\rho(p)$ the linear Hodge decomposition in $L^2$ of $A\nabla u$. Namely, there exist two functions $C$ and $D$, unique up to additive constants, elements of $W^{1,2}_0(B_\rho(p))$ and $W^{1,2}(B_\rho(p))$ respectively, and such that

\beq
\label{IV.14}
A\,\nabla u=\nabla C+\nabla^\perp D\quad .
\eeq
with
\be
\Delta C=div(A\,\nabla u) =\nabla^\bot B \cdot\nabla u
\ee
and 
\be
\Delta D=-\nabla A\cdot\nabla ^\perp u \quad.
\ee
Wente's lemma \ref{wente} guarantees that $C$ lies in $W^{1,2}$, and moreover
\beq
\label{IV.16}
\int_{B_\rho(p)}|\nabla C|^2 \, dz \le C\left( \int_{B_\rho(p)}|\nabla B|^2\, dz \right)\left( \int_{B_\rho(p)}|\nabla u|^2 \, dz\right)\quad.
\eeq 
Then, we introduce the decomposition
$D=\phi+v$, with $\phi$  satisfying
\beq
\label{IV.18}
\lf\{
\begin{array}{l}
\ds\Delta\phi=-\nabla A\cdot\nabla ^\perp u\quad\quad\mbox{ in }\quad B_\rho(p)\\[5mm]
\ds \quad\phi=0\quad\quad\quad\mbox{ on }\quad\p B_\rho(p)\quad ,
\end{array}
\rg.
\eeq
and with $v$ being harmonic. Once again, Wente's lemma \ref{wente} gives us the estimate
\be
\label{IV.19}
\int_{B_\rho(p)}|\nabla\phi|^2\, dz \le C\left(\int_{B_\rho(p)}|\nabla A|^2\, dz \right)\left( \int_{B_\rho(p)}|\nabla u|^2\, dz \right)\quad.
\ee
Using the fact that $\ds \rho\mapsto \frac{1}{\rho^2}\int_{B_\rho(p)} \vert \nabla v\vert^2 \, dz$ for any harmonic function, see lemma II.1 of \cite{Ri6}. We get, for any $0\leq \delta \leq1$,  that
$$
\int_{B_{\delta \rho}(p)} \vert \nabla v\vert^2 \, dz \leq \delta^2 \int_{B_{\rho}(p)} \vert \nabla v\vert^2 \, dz .$$
Finally, we have  

\beq
\label{IV.20}
\begin{array}{rl}
\ds \int_{B_{\delta \rho}(p)}|\nabla D|^2\, dz &\ds\le 2\delta^2\int_{B_\rho(p)}|\nabla D|^2\, dz\\[5mm]
 &\ds\ +2\int_{B_\rho(p)}|\nabla\phi|^2\, dz\quad .
\end{array}
\eeq
Bringing altogether (\ref{IV.14}), (\ref{IV.16}), and (\ref{IV.20}) produces
\beq
\label{IV.21}
\begin{array}{rl}
\ds\int_{B_{\delta \rho}(p)}|A\,\nabla u|^2\, dz&\ds\le2\delta^2\int_{B_\rho(p)}|A\,\nabla u|^2\, dz\\[5mm]
 &\ds \quad+C\,\eps_0\ \int_{B_\rho(p)}|\nabla u|^2\, dz .
 \end{array}
 \eeq
Using the hypotheses that $A$ and $A^{-1}$ are bounded in $L^\infty$, it follows from (\ref{IV.21}) that for all  $0<\delta<1$, there holds the estimate
\beq
\label{IV.22}
\begin{array}{rl}
\ds\int_{B_{\delta \rho}(p)}|\nabla u|^2\, dz&\ds\le2\|A^{-1}\|_\infty \,\|A\|_\infty\delta^2\int_{B_\rho(p)}|\nabla u|^2\, dz\\[5mm]
 &\ds \quad+C\,\|A^{-1}\|_\infty\eps_0\ \int_{B_\rho(p)}|\nabla u|^2\, dz\quad.
 \end{array}
 \eeq
Next, we choose $\eps_0$ and $\delta$ strictly positive, independent of $\rho$ et $p$, and such that
\[
2\|A^{-1}\|_\infty \,\|A\|_\infty\delta^2+C\,\|A^{-1}\|_\infty\eps_0=\frac{1}{2}\quad.
\]
For this particular choice of $\delta$, we have thus obtained the inequality
\[
\int_{B_{\delta \rho}(p)}|\nabla u|^2\, dz\le \frac{1}{2}\int_{B_\rho(p)}|\nabla u|^2\, dz\quad.
\]
Classical results then yield the existence of some constant $\al>0$ for which
\[
\sup_{p\in B_{1/2}(0)\;,\; 0<\rho<\frac{1}{2}}\rho^{-\al}\int_{B_\rho(p)}|\nabla u|^2\, dz<+\infty\quad,
\]
which prove the $\eps$-regularity as already remark above.\\

Then, the energy gap follows easily remarking that, thanks to the conformal invariance, for all $R>0$ and some $q>2$, we have
$$\Vert  \nabla u\Vert_{L^q(B_{\frac{R}{4}})} \leq \frac{C_q}{R^\frac{q-2}{q}} \Vert \nabla u \Vert_{L^2(B_R)}.$$
Finally, the weak limit property is a just a special case of theorem I.5 of \cite{Ri3} which is one of the many consequences of theorem \ref{lx}.
\hfill $\Box$\\

Then we are in position to prove theorem~\ref{lp} which is the main result of this section once we will have established the following lemma.

\begin{lemma}
\label{L2infini}
There exists  $\delta>0$ such that for all  $r,R\in\R^*_+$ satisfying $2r<R$ for all $\Omega\in L^2(B_R\setminus B_r, so(n)\otimes\R^n)$ and $u\in W^{1,2}(B_R\setminus B_r, \R^n)$ satisfying
$$-\Delta u =\Omega\cdot \nabla u$$
and
$$ \sup_{r<\rho<\frac{R}{2}} \int_{B_{2\rho}\setminus B_\rho}\vert\Omega\vert^2 \, dz \leq \delta .$$
Then  there exists $C>0$, independent of $u$, $r$ and $R$, such that 

\beq
\label{i0}
\left\Vert \nabla u\right\Vert_{L^{2,\infty}(B_{R}\setminus B_r)}\leq  C \ \left[\sup_{ r<\rho<\frac{R}{2}} \int_{B_{2\rho}\setminus B_\rho}\vert \nabla u\vert^2 \, dx\right]^{1/2} .
\eeq
\end{lemma}

\medskip

\noindent {\it Proof of lemma~\ref{L2infini}.}\\

Let 
\[
\eps:=\sup_{ r<\rho<\frac{R}{2}} \int_{B_{2\rho}\setminus B_\rho}\vert \nabla u\vert^2 \, dz
\]
We assume $\delta$ to be smaller than $\eps_0$ in the epsilon regularity result theorem~\ref{ereg} in such a way that for any $2r<\rho<R/4$ one has
\beq
\label{i1}
\left[\frac{1}{\rho^2}\int_{B_{2\rho}\setminus B_\rho}|\nabla u|^4(x)\ dx\right]^{\frac{1}{4}}\le C\ \frac{\sqrt{\eps}}{\rho}
\eeq
Let $\lambda>0$. Let  $f(x):=|\nabla u|$ in $B_{\frac{R}{2}}\setminus B_{2r}$ and $f=0$ otherwise, we have that 
\beq
\label{i2}
\forall \rho>0  \quad\quad\int_{B_{2\rho}\setminus B_\rho}\ f^4(x)\ dx\le C\ \frac{\eps^2}{\rho^2}
\eeq
For any $\rho>0 $ denote  
$$
U(\lambda,\rho):=\{x\in B_{2\rho}\setminus B_\rho\ ;\ f(x)>\lambda\}\quad.
$$
With this notation, (\ref{i2}) implies that
\beq
\label{ii2}
\lambda ^4\ |U(\lambda,\rho)|\le C\ \frac{\eps^2}{\rho^2}.
\eeq
Let $k\in {Z}$ and $j\ge k$ we apply (\ref{ii2}) for $\rho:=2^{j}\lambda^{-1}$ and by summing over $j\ge k$ one obtains
\[
\lambda^2\ \left|\left\{x\in \R^2 \setminus B_{2^k\lambda^{-1}} \ ;\ f(x)>\lambda\right\}\right|\le C\ \sum_{j=k}^\infty2^{-2j}\ \eps^2\le C\ 2^{-2k}\ \eps^2
\]
So we deduce that for any $k\in{\mathbb Z}$
\beq
\label{i3}
\lambda^2\ \left|\left\{x\in \R^2 \ ;\ f(x)>\lambda\right\}\right|\le C\ 2^{-2k}\ \eps^2+ \pi\ 2^{2k}
\eeq
Taking $2^{2k}\simeq\eps$ we obtain
\beq
\label{i4}
\left\Vert \nabla u\right\Vert_{L^{2,\infty}\left(B_{\frac{R}{2}}\setminus B_{2r}\right)}\leq  C \ \left[\sup_{ r<\rho<\frac{R}{2}} \int_{B_{2\rho}\setminus B_\rho}\vert \nabla u\vert^2 \, dx\right]^{1/2} .
\eeq
using now the triangular inequality for the norm $L^{2,\infty}$ and the fact that the $L^{2,\infty}$ norm of $\nabla u$ is controled by the $L^2$ norm of $\nabla u$ over respectively $B_R\setminus B_{\frac{R}{2}}$ and $B_{2r}\setminus B_r$, (\ref{i4}) implies (\ref{i0}) and lemma~\ref{L2infini} is proved.\hfill $\Box$

\medskip

\noindent{\it Proof of theorem~\ref{lp}:}\\

Let  $\eps_0>0$ be the one of the theorem \ref{lx}.\\

{\bf Step 1: We reduce the problem to an $L^{2,1}$ estimate}\\

Indeed, we use the duality $L^{2,1}-L^{2,\infty}$ in order to infer
\[
\int_{B_{\frac{R}{2}}\setminus B_{2r}}\left|\frac{1}{\rho}\frac{\partial u}{\partial \theta}\right|^2\ dx\le \left\|\frac{1}{\rho}\frac{\partial u}{\partial \theta}\right\|_{L^{2,1}\left(B_{\frac{R}{2}}\setminus B_{2r}\right)}\ 
\left\|\frac{1}{\rho}\frac{\partial u}{\partial \theta}\right\|_{L^{2,\infty}\left(B_{\frac{R}{2}}\setminus B_{2r}\right)}
\]
Combining this inequality with (\ref{i0}) we obtain
\beq
\label{i5}
\int_{B_{\frac{R}{2}}\setminus B_{2\rho}}\left|\frac{1}{\rho}\frac{\partial u}{\partial \theta}\right|^2\ dx\le C\ \left\|\frac{1}{\rho}\frac{\partial u}{\partial \theta}\right\|_{L^{2,1}\left(B_{\frac{R}{2}}\setminus B_{2r}\right)}
\left[\sup_{ r<\rho<\frac{R}{2}} \int_{B_{2\rho}\setminus B_\rho}\vert \nabla u\vert^2 \, dx\right]^{1/2} 
\eeq
Hence, thanks to duality, it suffices to control the $L^{2,1}$-norm of $ \frac{1}{\rho} \frac{\p u}{\p \theta}$ by the $L^2$ norm of $\nabla u$ in the annulus in order to prove the theorem.\\

{\bf Step 2: We prove the theorem assuming that 
$$\int_{B_R\setminus B_r}\vert \Omega\vert^2\, dz < \eps_0.$$}
We start by extending $\Omega$, setting
$$
\widetilde{\Omega}=\left\{\begin{array}{l}\Omega \hbox{ on¬† } B_R\setminus B_r \\[5mm]
0 \hbox{ on ¬†} B_r \quad.\end{array}\right.
$$
Hence, thanks to theorem \ref{lx},  there exists $\widetilde{A}\in W^{1,2}(B_R, Gl_n(\R))\cap L^\infty(B_R, Gl_n(\R))$ such that 
$$div(\nabla \widetilde{A} -\widetilde{A}\,\widetilde{\Omega})=0$$
and
\beq
\label{est1}
 \int_{B_R} (\vert \nabla \widetilde{A}\vert^2 +\vert \nabla \widetilde{A}^{-1}\vert^2)\, dz + dist(\{\widetilde{A},\widetilde{A}^{-1},\},SO(n))\leq C \int_{B_R} \vert \widetilde{\Omega} \vert^2 \, dz .
 \eeq
Then, thanks to Poincar\'e's lemma, there exists $\widetilde{B}\in  W^{1,2}(B_R(0), M_n(\R))$ such that
\beq
\label{A1}
\nabla \widetilde{A} -\widetilde{A}\,\widetilde{\Omega}=\nabla^\bot \widetilde{B}
\eeq
and, thanks to (\ref{est1}) and  (\ref{A1}), we get  
$$\Vert \nabla \widetilde{¬†B} \Vert_{L^2(B_R)}\leq C \Vert \Omega\Vert_{L^2(B_R\setminus B_r)},$$
here $C$ is a constant independent of $\Omega$. Hence, u satisfies
$$-div(\widetilde{A}\,\nabla u)= \nabla^\bot \widetilde{B} \cdot\nabla u \hbox{  on  }  B_R\setminus B_r .$$
We extend $u$ to $B_R$ by $\widetilde{u}$ using the Whitney's extension theorem, see \cite{AdFo} or \cite{Ste} for instance, then we get $\widetilde{u} \in W^{1,2}(B_R)$ such that 
\beq
\label{wint}
\int_{B_R} \vert \nabla \widetilde{u}\vert^2\, dz \leq C  \int_{B_R\setminus B_r} \vert \nabla u\vert^2\, dz .
\eeq
We consider the Hodge decomposition of $\widetilde{A}\nabla \widetilde{u}$ on $B_R$, i.e there exists  $C\in W^{1,2}_0(B_R)$ and  $D\in W^{1,2}(B_R)$ such that
\beq
\label{dec}
\widetilde{A}\,\nabla \widetilde{u} = \nabla C + \nabla^\bot D.
\eeq
 Moreover, thanks to  (\ref{wint}), we get 
$$\int_{B_R} \vert \nabla C\vert^2\, dz + \int_{B_R} \vert \nabla D\vert^2\, dz = \int_{B_R}\vert \widetilde{A}\, \nabla \widetilde{u} \vert^2 \, dz \leq C \ \int_{B_R\setminus B_r} \vert \nabla u\vert^2\, dz \quad. $$
Here we use the fact that $C$ vanishes on the boundary to get that 
$$\int_{B_R} \nabla C\cdot\nabla^\bot D \, dz =0.$$
Then, on $B_R\setminus B_r$ , $C$ satisfies
$$\Delta C=\nabla^\bot \widetilde{B}\cdot \nabla u.$$
As usual, we split as follows $C=v+\phi$ where $\phi \in W^{1,2}_0(B_R\setminus B_r)$ and $v \in W^{1,2}(B_R\setminus B_r)$ which satisfy
$$\Delta \phi =  \nabla^\bot \widetilde{B}\cdot \nabla u $$
and 
$$  \Delta v=0.$$
On the one hand, thanks to lemma \ref{l4} we get, for $0<\lambda<1$, that
$$
\Vert \nabla \phi \Vert_{L^{2,1}(B_R\setminus B_{\lambda^{-1}r})}\leq C(\lambda)\ \Vert \nabla \widetilde{B}\Vert_2\ \Vert \nabla u\Vert_2\quad .
$$
On the other hand, we decompose $v$ as a Fourier series,
$$v= c_0 + d_0 \log(\rho) + \sum_{n\in \Z^*} (c_n \rho^n +d_n \rho^{-n})e^{in\theta}.$$
since $\frac{1}{\rho} \frac{\p v}{\p \theta}$ has no logarithm part, we conclude as in lemma \ref{l3} that for any $0<\lambda<1$ we have

$$
\left\Vert\frac{1}{\rho} \frac{\p v}{\p \theta}\right\Vert_{L^{2,1}(B_{\lambda R}\setminus B_{\lambda^{-1}r})}\leq C(\lambda)\ \Vert \nabla v \Vert_2 .
$$
The Dirichlet principle implies that 
$$\Vert \nabla v \Vert_2 \leq \Vert \nabla C \Vert_2,$$
then we get
\beq
\label{C}
\left\Vert\frac{1}{\rho} \frac{\p C}{\p \theta}\right\Vert_{L^{2,1}(B_{\lambda R}\setminus B_{\lambda^{-1}r})}\leq C(\lambda)\ \Vert \nabla u \Vert_{L^2(B_R\setminus B_r)} .
\eeq
Now we estimate $D$, which satisfies the following equation 
$$\Delta D =\nabla \widetilde{A}\cdot \nabla^\bot \widetilde{u} \hbox{   on   }B_R.$$
Then, we also decompose $D$ as $D=v+\phi$ where $\phi \in W^{1,2}_0(B_R)$ and $v \in W^{1,2}(B_R)$. 
$$\Delta \phi =  \nabla \widetilde{A}\cdot \nabla^\bot u $$
and 
$$  \Delta v=0.$$
In the one hand, thanks to lemma \ref{wente}, we have
\be
\begin{split}
\Vert \nabla \phi\Vert_2 &\leq \Vert \nabla \phi\Vert_{L^{2,1}(B_R)} \leq C\ \Vert \nabla \widetilde{A}\Vert_2\ \Vert \nabla \widetilde{u}\Vert_2\\[5mm]
&\leq C\ \Vert \nabla u \Vert_{L^2(B_R\setminus B_r)} .
\end{split}
\ee
in the other hand, since $v$ is harmonic, for any $0<\lambda<1$, we have
\be
\begin{split}
\Vert \nabla v\Vert_{L^{2,1}(B_{\lambda R})} &\leq C(\lambda)\  \Vert \nabla v\Vert_{L^2(B_{ R})} \\[5mm]
&\leq C(\lambda)\  \Vert \nabla D\Vert_{L^2(B_{ R})} \\[5mm]
 &\leq C(\lambda)\  \Vert \nabla u\Vert_2 .
\end{split}
\ee
Finally
\beq
\label{D}
\Vert \nabla D \Vert_{L^{2,1}(B_{\lambda R}\setminus B_{\lambda^{-1}r})} \leq  C(\lambda)\ \Vert \nabla u\Vert_2 .
\eeq
Combining (\ref{dec}), (\ref{C}) and (\ref{D}), we get 
$$\left\Vert \widetilde{A}\, \frac{1}{r} \frac{\p \widetilde{u}}{\p \theta} \right\Vert_{L^{2,1}(B_{\lambda R}\setminus B_{\lambda^{-1}r})} \leq  C(\lambda)\  \Vert \nabla u\Vert_2.$$
Finally, using (\ref{est1}), we get that 
\beq
\label{fs2}
\left\Vert \frac{1}{\rho} \frac{\p \widetilde{u}}{\p \theta} \right\Vert_{L^{2,1}(B_{\lambda R}\setminus B_{\lambda^{-1}r})} \leq  C(\lambda)\  \Vert \nabla u\Vert_2,
\eeq
which proves,  as remark at the hand of step 1, the theorem under the extra assumption.\\

{\bf Step 3: General case}\\

We construct two sequences of radii $r_i$ and $R_i$ such that
$$r=r_0<r_1=R_0<\dots < r_{i+1}=R_i <\dots< R_N=R$$ 
with
$$\int_{B_{R_i}\setminus B_{r_i}} \vert \Omega\vert^2 \, dz \leq \eps_0  $$
and
$$N\leq \frac{\int_{B_{R}\setminus B_{r}} \vert \Omega\vert^2 \, dz}{\eps_0}  .$$
First, applying (\ref{fs2}) of step 2, we get that 
\beq
\label{f1}
 \left\Vert  \frac{1}{\rho} \frac{\p u}{\p \theta} \right\Vert_{L^{2,1}(B_{\lambda R_i}\setminus B_{\lambda^{-1} r_i})} \leq  C(\lambda)\ \Vert \nabla u\Vert_{L^2(B_{R_i}\setminus B_{r_i})}
 \eeq
We chose $\delta$ such that 
$$\delta< \frac{\eps_0}{4}$$
hence for all $i$ we have
$$\int_{B_{4r_i}\setminus B_{\frac{r_i}{4}}} \vert \Omega\vert^2\, dz <4\delta<\eps_0$$
Let $S_i=\min(R,4r_i)$ and $ s_i =\max(r,\frac{r_i}{4})$, then we apply again (\ref{fs2}) of step 2 on $B_{S_i}\setminus B_{s_i}$, which gives 
\beq
\label{f2}
\left\Vert  \frac{1}{\rho} \frac{\p u}{\p \theta} \right\Vert_{L^{2,1}(B_{\lambda S_i}\setminus B_{\lambda^{-1} s_i})} \leq  C(\lambda)\ \Vert \nabla u\Vert_{L^2(B_{S_i}\setminus B_{s_i})}
\eeq
Finally, summing (\ref{f1})  and (\ref{f2}),  for $i=0$ to $N$, we get 
$$ \left\Vert  \frac{1}{\rho} \frac{\p u}{\p \theta} \right\Vert_{L^{2,1}(B_{\lambda R}\setminus B_{\lambda^{-1}r})} \leq  C(\lambda)\ \Vert \nabla u\Vert_2,$$
which achieves the proof of  theorem~\ref{lp}      .\hfill$\square$\\

We shall now make use of the theorem~\ref{lp} in order to prove the quantization of the angular part of the energy for solutions to antisymmetric elliptic systems.\\

We wil call a {\bf bubble} a solution $u\in W^{2,1}(\R^2,\R^n)$ of the equation
$$-\Delta u = \Omega \cdot \nabla u \hbox{ on } \R^2,$$
where $\Om\in L^2(\R^2,so(n)\otimes{\mathbb R}^2)$.\\

{\it Proof of theorem~\ref{th-I.1} :}\\

First we are going to separate $B_1$ in three parts: one where $u_k$ converge to a limit solution, some neighborhoods where the energy concentrates and where blow some bubbles and some neck regions which join the first  two parts. This ''bubble-tree'' decomposition is by now classical, see \cite{Pa} for instance, hence we just sketch briefly how to proceed.\\ 
 
{\bf Step 1 : Find the point of concentration}\\

Let  $\eps_0$ be  the one of theorem \ref{ereg} and $\delta$ the one of theorem \ref{lp}. Then, thanks to (\ref{be}), we easily proved that there exists finitely many points $a^1,\dots, a^n$ where 
\beq
\label{cp}
\int_{B(a_i, r)}  \vert \Omega_k\vert^2 \, dz \geq \eps_0 \hbox{ for all } r>0.
\eeq
Moreover, using theorem \ref{ereg}, we prove that  there exists $\Om_\infty\in L^2(B_1,so(n)\otimes{\mathbb R}^2)$  and $u_\infty\in W^{2, 1}(B_1,\R^n)$ a solution of $-\Delta u=\Omega_\infty \cdot \nabla u$ on $B_1$, such that, up to a subsequence,  
$$\Omega_k \rightharpoonup \Omega_\infty \hbox{ in }  L^2_{loc}(B_1,so(n) \otimes \R^2),$$
and
$$u_k\rightarrow u_\infty \hbox{ on } W^{1,p}_{loc}(B_1\setminus \{ a^1,\dots,a^n\}) \hbox{ for all } p\geq1 .$$
Of course, if $\Vert \Omega_k\Vert_\infty =O(1)$ or  $\Omega_k =O(\nabla u_k)$, then $u_k$ is bounded in $W^{2,\infty}$ which gives the convergence in  $C^{1,\eta}_{loc}$ for all $\eta\in[0,1[$.\\ 

{\bf Step 2 : Blow-up around $a^i$ }\\

We choose $r_i>0$ such that 
$$\int_{B(a^i, r^i)} \vert \Omega_\infty \vert^2 \, dz \leq \frac{\eps_0}{4} .$$  
Then, we define a center of mass of  $B(a^i, r^i)$ with respect to $\Omega_k$ in the following way
$$a_k^i= \left(\frac{\displaystyle\int_{B(a^i,r^i)} x^\alpha \vert \Omega_k\vert^2\, dz}{\displaystyle\int_{B(a^i,r^i)}\vert \Omega_k\vert^2\, dz} \right)_{\alpha=1,2}.$$
Let $\lambda^i_k$ be a positive real such that 
 $$\int_{B(a^i_k, r^i)\setminus B(a^i_k, \lambda_k^i)} \vert \Omega_k \vert^2 \, dz = \min\left(\delta, \frac{\eps_0}{2}\right).$$  
 We set $\widetilde{u}_k(z) = u_k( a^i_k+ \lambda_k^i z)$, $\widetilde{\Omega}_k(z) = \lambda_k^i \,\Omega_k( a^i_k+ \lambda_k^i z)$ and $N^i_k =B(a_k^i,r^i)  \setminus B(a_k^i, \lambda_k^i)$.\\
Observe that the scaling we chose for defining $\widetilde{\Omega}_k(z)$ guaranties that
\[
\int_{B(0,r^i/\lambda_k^i)}|\widetilde{\Omega}_k|^2\ dz\le \int_{B(a_k^i,r^i)}|\Omega_k|^2\ dx\le C<+\infty
\]
moreover we have
\[
-\Delta\widetilde{u}_k=\widetilde{\Omega}_k\cdot \nabla \widetilde{u}_k\quad.
\]
Modulo extraction of a subsequence, we can assume that for each $i$
\[
\nabla\widetilde{u}_k\rightharpoonup \nabla\widetilde{u}_\infty\quad\quad\mbox{ weakly in }L^2_{loc}(\R^2,{\mathbb R}^n)\quad\quad\widetilde{\Omega}_k\rightharpoonup\widetilde{\Omega}_\infty\mbox{ weakly in }L^2_{loc}( \R^2,so(n)\otimes{\mathbb R}^2)
\]
The {\bf weak limit property} of theorem \ref{ereg} implies that $\widetilde{u}_\infty$ and $\widetilde{\Omega}_\infty$ satisfy what we call a {\it bubble } equation
\[
-\Delta \widetilde{u}_\infty=\widetilde{\Omega}_\infty \cdot\nabla \widetilde{u}_\infty\quad.
\]

{\bf Step 3 : Iteration}\\
 
 Two cases have to be considered separatly:\\

Either $\widetilde{\Omega}_k$ is subject to some concentration phenomena as (\ref{cp}), and then we find some new points of concentration and  a limiting solution which is  a bubble (possibly trivial), in other words a solution of  $-\Delta u=\widetilde{\Omega}_\infty \cdot\nabla u$, where $\widetilde{\Omega}_\infty$ is a weak limit of $\widetilde{\Omega}_k$. In such a case we apply step $2$ to our new concentration points.\\ 

Or, $\widetilde{u}_k$  converges uniformly on every compact subset to a bubble (possibly trivial).\\

Of course this process has to stop, since we are assuming a uniform bound on $\Vert \Omega_k\Vert_2$  and  each step is consuming at least $\min\left(\delta, \frac{\eps_0}{2}\right)$ of energy of $\Omega_k$. This process is sketched in the following picture.

\begin{figure}[!h]
\centering
\input{laurain-fig1.tex} 
\caption{Decomposition of $\Sigma$}
\label{parti}
\end{figure}
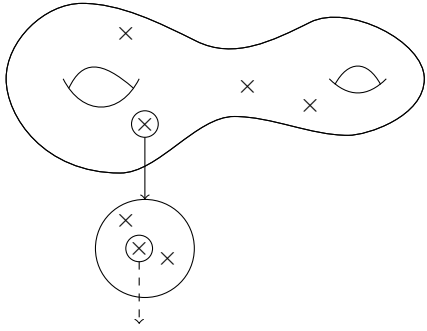

{\bf Analysis of a neck region: }\\

A neck region is an annullar region $N_k^i= B\left(a_k^i,\mu^i_k \right )\setminus B\left (a_k^i, \lambda_k^i\right)$ such that 
$$\lim_{k\rightarrow +\infty} \frac{\lambda_k^i} {\mu^i_k}=0,$$
\beq
\label{aa1}
\int_{N_k^i} \vert \Omega_k\vert^2 \, dz \leq \min\left(\delta, \frac{\eps_0}{2}\right)
\eeq
and
$$X_k= \nabla^\bot d(a_k^i,\, .\,).$$
In order to prove  theorem~\ref{th-I.1}, we start by proving  a weak estimate on the energy of gradient in the neck region . First we remark that, for all $\eps>0$, there exists $r>0$ such that for all $\rho>0$ such that 
$$B_{2\rho}(a_k^i)\setminus B_\rho(a_k^i) \subset N_k^i(r)$$
where $N_k^i(r)= B\left(a_k^i,r\mu^i_k \right )\setminus B\left (a_k^i, \frac{\lambda_k^i}{r}\right)$, we have 
\beq \label{we1}
\int_{B_{2\rho}(a_k^i)\setminus B_\rho(a_k^i) } \vert \nabla u\vert^2 \, dz \leq \eps\quad .
\eeq
If this would not be the case there would exist a sequence $\rho_k^i\rightarrow 0$ such that, up to a subsequence, $\hat{u}_k = u_k(a_k^i+ \rho_k^i z)$ converge 
with respect to every $W^{1,p}-$norm to a non-trivial solution of
$$ -\Delta \hat{u} = \widehat{\Omega}_\infty \cdot \nabla \hat{u} \hbox { on } \R^2\setminus\{ 0\}, $$
where $\widehat{\Omega}_\infty$ is a weak limit, up to a subsequence, of $\widehat{\Omega}_k$.
Using the fact that the $W^{1,2}$-norm of $\hat{u}_k$ is bounded, we deduce using Schwartz lemma that it has to be in fact a solution on the whole plane.
Using this time  the second part of theorem \ref{ereg} we deduce that $\hat{\Omega}_\infty$  have energy at least $\eps_0$, which contradicts (\ref{aa1}).\\

Finally, using theorem \ref{lp} on each $N_k^i(r)$ , we obtain
\be
\begin{split}
\lim_{r\rightarrow 0}\lim_{k\rightarrow +\infty}  \left\Vert \langle \nabla  u_k, X_k\rangle \right\Vert_{L^{2}(N_k^i(r))} &\leq C\ \lim_{r\rightarrow 0} \lim_{k\rightarrow +\infty}  \left(\sup_{\rho} \int_{B_{2\rho}(a_k^i)\setminus B_\rho(a_k^i) } \vert \nabla u\vert^2 \, dz\right)   \\
&=0.
\end{split}
\ee
Which achieves the proof of  theorem~\ref{th-I.1}.\hfill$\square$\\

This phenomena of quantization of the angular part of the gradient seems to be quite general for systems with antisymmetric potentials, in a forthcoming paper \cite{LaRi} we investigate the quantization for some fourth order elliptic systems in 4-dimension. 

\section{Energy Quantization for critical points to conformally invariant lagrangians.}
In the present section we are going to use theorem~\ref{th-I.1} in order to prove theorem~\ref{th-I.3}\\

In his proof of the Heinz-Hildebrandt's regularity conjecture, the second author prove that the Euler Lagrange equations to general conformally invariant lagrangians which are coercive
and of quadratic growth can be written in the form of an elliptic system with an antisymmetric potential.  Precisely we have
\begin{thm}[Theorem I.2 \cite{Ri3}] 
 Let $N^k$  be a $C^2$ submanifold of $\R^m$ and $\omega$ be a $C^1$ $2$-form on $N^k$ such that the $L^\infty$-norm of $d\omega$ is bounded on $N^k$. Then every critical point in $W^{1,2}(B_1,N^k)$ of the Lagrangian
\beq
\label{cil}
F(u) =\int_{B_1} \left[\vert \nabla u\vert^2 + u^\ast\omega\right]\, dz
\eeq
satisfies 
$$ -\Delta u = \Omega\cdot \nabla u,$$
with
\beq
\label{aa2}
\Omega^i_j=[A^i(u)_{j,l}-A^j(u)_{i,l}]\ \nabla u^l + \frac{1}{4} [H^i(u)_{j,l}-H^j(u)_{i,l}]\ \nabla^\bot  u^l\quad.
\eeq
where $A$ and $\lambda$ are in $C^0(B_1, M_m(\R)\otimes \bigwedge^1\R^2)$ satisfy
$$\sum_{j=1}^m A^j_{i,l} \nabla u^j =0$$
and $H^i_{j,l}:= d(\pi^\ast\omega)(\eps_i,\eps_j,\eps_l)$ where, in a neighborhood
of $N^k$, $\pi$ is the orthogonal projection onto $N^k$ and $(\eps_i)_{i=1\cdots m}$ is the canonical basis of ${\R}^m$.\hfill $\Box$
 \end{thm}

 From (\ref{aa2}) we observe that  for critical points to a conformally invariant $C^1-$Lagrangian, there exists 
  \beq 
 \label{H1}
 \Lambda \in C^0(TN\otimes \R^2, so(n)\otimes \R^2)
 \eeq
 such that 
 \beq
 \label{H2}
\Lambda (v) =O(\vert v\vert),
 \eeq
 moreover we remark that $\Lambda(u,\nabla u)\cdot \nabla u $ is always  orthogonal to $\nabla u $ in the following sense
 \beq
 \label{H3}
 \left\langle  \frac{\partial u}{\partial x_k} , \Lambda(u,\nabla u)\cdot \nabla u\right\rangle =0 \hbox{ for } k=1,2.
 \eeq
For $\Lambda \in C^0(TN\otimes \R^2, so(n)\otimes \R^2)$, we call a {\bf $\Lambda$-bubble}  a solution $\omega\in W^{2,1}(\R^2,\R^n)$ of the equation
$$-\Delta \omega = \Lambda(\omega, \nabla \omega)\cdot \nabla \omega \hbox{ on } \R^2.$$\\
 
\begin{thm}
\label{qw} Let $u_k\in W^{1,2}(B_1,\R^n)$ be a sequence of critical points of a functional which is conformally invariant, which satisfies
\beq
\label{eq2}
-\Delta u_k =\Lambda(u_k, \nabla u_k)\cdot \nabla u_k ,
\eeq 
where $\Lambda$ satisfies (\ref{H1}), (\ref{H2}) and (\ref{H3}). Moreover we assume that $u_k$ has a bounded energy, i.e.
\be
\Vert \nabla u_k\Vert_{2} \leq M.
\ee
Then there exists $u_\infty\in W^{1, 2}(B_1,\R^n)$ a solution of $-\Delta u_\infty=\Lambda(u_\infty,\nabla u_\infty)\cdot \nabla u_\infty$ on $B_1$, 
 $l\in\N^*$ and
\begin{enumerate}
\item $\omega^1,\dots, \omega^l$ some non-constant $\Lambda$-bubbles
\item $a_k^1,\dots,a_k^l$ a family of converging sequences of points of $B_1$
\item  $\lambda_k^1,\dots,\lambda_k^l$ a family of  sequences of positive reals converging all to zero.
\end{enumerate}
such that, up to a subsequence, 
$$u_k\rightarrow u^\infty \hbox{ on } C^{1,\eta}_{loc}(B_1 \setminus \{ a^1_\infty,\dots,a^l_\infty\}) \hbox{ for all }  \eta\in[0,1[$$
and 
$$ \left\Vert \nabla \left( u_k -u_\infty -\sum_{i=1}^l \omega^i_k \right) \right\Vert_{L^2_{loc}(B_1)}\rightarrow 0\quad,$$
where $\omega_k^i=\omega(a_k^i + \lambda_k^i z)$. \hfill$\Box$
\end{thm}

Since (\ref{H3}) holds for any system issued from a lagrangian of the form (\ref{cil}), it is clear that theorem~\ref{th-I.3} is a consequence of
theorem~\ref{qw}.

\medskip

\noindent{\it Proof of theorem~\ref{qw}    :}\\

From the previous section, we have the quantization of the angular part of the gradient. To prove theorem~\ref{qw} it suffices then to prove the energy quantization for the radial part 
of the energy. Since $u_k$ satisfies (\ref{eq2}) then $u_k\in W^{2,p}(B_{\mu^i_k}(a_k^i))$ for all $p<\infty$, see theorem IV.3 of \cite{Ri6} or lemma 7.1 of \cite{ShTo}, hence we can multiply (\ref{eq2}) by $\rho \frac{\p u_k}{\p \rho}$ and integrate. Using  (\ref{H3}) we have, for any $r\in [0, \mu^i_k]$, 

$$0  = \int_{B_r}\langle  \rho \frac{\p u_k}{\p \rho} , \Omega\cdot \nabla u_k\rangle \, dz  = \int_{B_r} \langle \rho \frac{\p u_k}{\p \rho}, \Delta u_k\rangle\, dz. $$
Using Pohozaev identity, we get for all $r \in[0, \mu^i_k]$
 $$\int_{\p B_r} \left\vert \frac{\p u_k}{\p \rho}\right\vert^2 \, d\sigma =  \int_{\p B_r} \left\vert \frac{1}{\rho} \frac{\p u_k}{\p \theta}\right\vert^2 \, d\sigma$$
Finally, we have 
$$\lim_{r\rightarrow 0}\lim_{k\rightarrow +\infty}  \Vert  \nabla u_k \Vert_{L^{2}(N_k^i(r))} =0,$$
which concludes the proof of the theorem. \hfill$\square$\\

In particular we get the quantization for the solution of the problem of prescribed mean curvature. Indeed, an immersion of a Riemann surface $\Sigma$ into $\R^3$ with prescribed mean curvature $H\in C^0(\R^3,\R)$ satisfies the following $H$-system 
\beq
\label{*}
\Delta u= 2H(u)\, u_x \wedge u_y,
\eeq
where $z=x+iy$ are some local conformal coordinates on $\Sigma$.\\ 

In order to state precisely our theorem, we define the notion of $H$-bubble as being a map  $\omega \in W^{1,2}(\R^2,\R^3)$ satisfying 
$$\Delta \omega= 2H(\omega)\,  \omega_x \wedge \omega_y  \hbox{ on } \R^2.$$
We shall also rescale the Riemann surface around a point. To that aim we will introduce some conformal chart. Precisely  there exists $\delta>0$ such that  for any $a\in \Sigma$ and $0<\lambda<\delta$  there exists a map $\Phi_{a,\lambda}:B(a,\delta)\rightarrow \R^2$ which is a conformal-diffeomorphism, sends $a$ to $0$ and $B(a,\lambda)$ to $B(0,1)$. We also associate to each point a cut-off function  $\chi_{a} \in C^{\infty}(\Sigma)$ which satisfies

$$
\left\{\begin{array}{l} \chi_{a} \equiv 1 \hbox{ on } B(a,\frac{\delta}{2}) \\[5mm] \chi_{a} \equiv 0 \hbox{ on } \Sigma\setminus B(a,\delta)\quad.  \end{array}\right. $$ 

\begin{cor}
\label{qh} Let $\Sigma$ be a closed Riemann surface, $H\in C^0 (\R^3,\R)$ and $u_k\in W^{2, 1}(\Sigma,\R^3)$ a sequence of non-constant solution of (\ref{*}) on $\Sigma$ then there exists,
$u_\infty\in W^{2, 1}(\Sigma,\R^3)$ a solution of (\ref{*}) , 
 $k\in \N^*$ and
\begin{enumerate}
\item $\omega^1,\dots, \omega^l$ a family of $H$-bubbles
\item $a_k^1,\dots,a_k^l$ a family of converging sequences of point of $\Sigma$
\item  $\lambda_k^1,\dots,\lambda_k^l$ a family of  sequences of positive reals converging all to zero
\end{enumerate}

$$u_k\rightarrow u^\infty \hbox{ on } C^{1,\eta}_{loc}(\Sigma\setminus \{ a_1^\infty,\dots,a_k^\infty\}) \hbox{ for all } \eta\in [0,1[$$
and moreover
$$ \left\Vert \nabla \left( u_k -u_\infty -\sum_{i=1}^l \chi_{a^i_k} \left(\omega^i \circ \Phi_{a^i_k,\lambda^i_k}\right)\right) \right\Vert_2\rightarrow 0\quad.$$
\end{cor}

\medskip

We end up this section by mentioning a recent work by Da Lio, \cite{DaL} in which energy quantization results for fractional harmonic maps (which are also conformally invariant
in some dimension) are established using also Lorentz space uniform estimates.

\section{Other applications to pseudo-holomorphic curves, harmonic maps and Willmore surfaces}

In this section we give some more applications of the uniform Lorentz-Wente estimates of section 2 to problems where the conformal invariance play again a central role.\\

In the present section we are interested  with Wente's type estimate for first order system of the form
\beq
\label{5.1}
 \nabla \phi =\sum_{i=1}^n a_i \,\nabla^\bot b_i 
 \eeq
 Taking the divergence of this system gives the classical order 2 Wente system
 \beq
 \label{5.2}
  \Delta \phi =\sum_{i=1}^n \nabla a_i \cdot\nabla^\bot b_i 
 \eeq 
 The gain of information provided by a first order system of the form (\ref{5.1}) in comparison to
 classical second order system (\ref{5.2}) is illustrated by the fact that, in the first order case, no assumption on the behavior of the solution $\phi$ at the boundary of the annulus is needed
 in order to obtain the Lorentz - Wente type estimates of section 2. This is proved in lemma~\ref{fo}. This fact can be applied to geometrically interesting situations that we will
 describe at the end of the present section.
\subsection{Lorentz-Wente type estimates for first order Wente type equations.}

The goal of this subsection is to prove the following lemma.

\begin{lemma}
\label{fo}
Let $n\in\N^*$, $(a_i)_{1\leq i\leq n}$ and $(b_i)_{1\leq i\leq n}$ be two families of maps in $W^{1,2}(B_1)$, $0<\eps<\frac{1}{4}$ and $\phi\in W^{1,2}(B_1\setminus B_\eps)$ which satisfies
\beq
\label{z1}
 \nabla \phi =\sum_{i=1}^n a_i \,\nabla^\bot b_i .
 \eeq
 Then, for $0< \lambda<1$, there exists a positive constant  $C(\lambda)$ independent of $\phi, a_i$ and $b_i$ such that 
 $$\Vert \nabla \phi\Vert_{L^{2,1}( B_{\lambda}\setminus  B_{\lambda^{-1}\eps })} \leq C(\lambda) \left(  \sum_{i=1}^n \Vert\nabla  a_i\Vert_2\  \Vert \nabla  b_i\Vert_2 + \Vert \nabla \phi\Vert_2\right) .$$
 \hfill $\Box$
\end{lemma}

\noindent{\it Proof of lemma~\ref{fo} :}\\

Taking the divergence of (\ref{z1}), gives
$$ 
\Delta \phi =\sum_{i=1}^n \nabla a_i\cdot \nabla^\bot b_i \quad.
$$
Hence, as in the previous lemma, we start by considering a solution of this equation on the whole disk and equal to zero on the boundary. Let $\varphi\in W^{1,1}_0(B_1)$ be the solution of
$$ 
\Delta \varphi = \sum_{i=1}^n \nabla a_i \cdot\nabla^\bot b_i \quad.
$$
Then, thanks to the improved Wente's inequality (\ref{w1}), we have 
\beq
\label{y3}
\Vert\nabla \varphi\Vert_{L^{2,1}(B_1)} \leq C\ \sum_{i=1}^n \Vert \nabla a_i\Vert_{2}\  \Vert \nabla b_i\Vert_2 \quad.
\eeq
We now consider the difference $v=\phi-\varphi$, which is an harmonic function on $B_1\setminus B_\eps$. Following the proof of the lemma \ref{l3}, it suffices to control  the logarithmic part of the decomposition in Fourier series. To that aim we set 
$$\overline{\phi}(\rho) =  \frac{1}{2\pi} \int_0^{2\pi} \phi(\rho,\theta)\, d\theta . $$ 
We have
\be
\begin{split}
\frac{d{\overline{\phi}}}{d\rho}&= \frac{1}{2\pi} \int_0^{2\pi} \frac{\partial \phi}{\partial \rho}(\rho,\theta)\, d\theta=\frac{1}{2\pi} \sum_{i=1}^n  \int_0^{2\pi} a_i \frac{\partial b_i }{\partial \theta} \, \frac{d\theta}{\rho} \\
&=\frac{1}{2\pi} \sum_{i=1}^n  \int_0^{2\pi} (a_i-\overline{a_i}) \frac{\partial b_i }{\partial \theta} \, \frac{d\theta}{\rho} 
\end{split}
\ee
Hence
\be
\begin{split}
\left\vert \frac{d{\overline{\phi}}}{d\rho} \right\vert\leq \frac{1}{2\pi} \sum_{i=1}  \left(\int_0^{2\pi} \vert a_i-\overline{a_i}\vert^2  \, d\theta\right)^\frac{1}{2} \left(\int_0^{2\pi} \left\vert \frac{1}{\rho}\frac{\partial b_i }{\partial \theta}\right\vert^2 \, d\theta\right)^{\frac{1}{2}} .
\end{split}
\ee
Which gives, thanks to Poincar\'e's inequality on the circle,  
\be
\left\vert \frac{d{\overline{\phi}}}{d\rho} \right\vert \leq C\ \sum_{i=1}  \left(\int_0^{2\pi}\left\vert\frac{\partial a_i }{\partial \theta}\right\vert^2  \, d\theta\right)^\frac{1}{2} \left(\int_0^{2\pi} \left\vert \frac{1}{\rho}\frac{\partial b_i }{\partial \theta}\right\vert^2 \, d\theta\right)^{\frac{1}{2}}, 
\ee
where $C$ is a constant independent of $\phi$.\\

Then integrating over $[1,\eps]$, we get
\beq
\label{5a3}
\begin{split}
\int_\eps^1  \left\vert \frac{d{\overline{\phi}}}{d\rho} \right\vert \, d\rho 
&\leq C \sum_{i=1}^n \int_\eps^1    \left(\int_0^{2\pi} \left\vert\frac{\partial a_i }{\partial \theta}\right\vert^2  \, d\theta\right)^\frac{1}{2} \left(\int_0^{2\pi} \left\vert \frac{1}{\rho}\frac{\partial b_i }{\partial \theta}\right\vert^2 \, d\theta\right)^{\frac{1}{2}} d\rho\\
&\leq C \sum_{i=1}^n \left( \int_{\mathrm{D}(0,1)\setminus B_\eps}   \left\vert\frac{1}{\rho}\frac{\partial a_i }{\partial \theta}\right\vert^2  \, \rho\ d\rho\, d\theta\right)^\frac{1}{2} \left( \int_{\mathrm{D}(0,1)\setminus B_\eps}    \left\vert \frac{1}{\rho}\frac{\partial b_i }{\partial \theta}\right\vert^2 \, \rho\ d\rho\, d\theta\right)^{\frac{1}{2}} \\
&\leq C\left( \sum_{i=1}^n \Vert \nabla a_i\Vert_2\  \Vert \nabla b_i\Vert_2\right) .
\end{split}
\eeq  
Moreover, by duality, we obtain
\beq
\label{5a4}
\int_\eps^1  \left\vert \frac{d{\overline{\varphi}}}{d\rho} \right\vert \, d\rho \leq \left\Vert \nabla \varphi  \frac{1}{\rho}\right\Vert_1\leq  \left\Vert \nabla \varphi \right\Vert_{L^{2,1}} \left\Vert \frac{1}{\rho} \right\Vert_{L^{2,\infty}} \leq  C  \Vert \nabla \varphi \Vert_{L^{2,1}} .
\eeq
The combination of (\ref{y3}), (\ref{5a3}) and (\ref{5a4}) gives then
\beq
\label{y4} 
\int_\eps^1  \left\vert \frac{d{\overline{v}}}{d\rho} \right\vert \, d\rho \leq C\left( \sum_{i=1}^n \Vert \nabla a_i\Vert_2  \Vert \nabla b_i\Vert_2\right) 
\eeq
Following the approaches we used in the proofs of the various lemma in section 2, we decompose $v$ as a Fourier series, which gives

$$v(\rho,\theta) = c_0 + d_0\ln(\rho) +\sum_{n\in \Z^{*}} (c_n \rho^n+d_n \rho^{-n}) e^{in\theta}.$$
We have
\[
\overline{v}(\rho)= c_0 + d_0\ln(\rho)
\]
Thanks to (\ref{y4}), we get that 
\beq
\label{y5}
|d_0|\ \log\frac{1}{\eps}  \leq C\left( \sum_{i=1}^n \Vert \nabla a_i\Vert_2\  \Vert \nabla b_i\Vert_2\right) \quad.
\eeq
We have moreover
\beq
\label{y5b}
\begin{split}
\|\nabla \overline{v}\|_{L^{2,1}(B_1\setminus B_\eps)}\simeq |d_0|\int_0^\infty|\{ x\in B_1\setminus B_\eps\ ;\ |x|^{-1}>t\}|^{1/2}\ dt\\[5mm]
=|d_0|\ \int_0^\infty|(B_1\setminus B_\eps)\cap B_{1/t}|^{1/2}\ dt\le \pi\, |d_0|\ \int_0^{1/\eps}\ \frac{dt}{\max\{t,1\}}=\pi\,|d_0|\ \left[1+\log\frac{1}{\eps}\right]
\end{split}
\eeq
Thus combining (\ref{y5}) and (\ref{y5b}) we have in one hand
\beq
\label{y5c}
\|\nabla \overline{v}\|_{L^{2,1}(B_1\setminus B_\eps)}\le C\left( \sum_{i=1}^n \Vert \nabla a_i\Vert_2\  \Vert \nabla b_i\Vert_2\right) \quad,
\eeq
in the other hand, as in lemma \ref{l3}, we have
\beq
\label{y6}
\left\Vert  \sum_{n\in \Z^{*}} (c_n\rho^n+d_n \rho^{-n}) e^{in\theta}   \right\Vert_{L^{2,1}( B_{\lambda}\setminus  B_{\lambda^{-1}\eps })}\leq C(\lambda) \Vert \nabla v\Vert_2 \leq C(\lambda)\ \Vert \nabla \phi \Vert_2  .
\eeq
Combining (\ref{y5c}), (\ref{y6}) we have for any $\lambda\in (0,1)$ the existence of a positive constant $C(\lambda)>0$ such that
\beq
\label{y7}
\Vert \nabla v\Vert_{L^{2,1}( B_{\lambda}\setminus  B_{\lambda^{-1}\eps })} \leq C(\lambda) \left(  \sum_{i=1}^n \Vert\nabla  a_i\Vert_2\  \Vert \nabla  b_i\Vert_2 + \Vert \nabla \phi\Vert_2\right)\ .
\eeq
Finally summing (\ref{y3}) and (\ref{y7}) gives the desired inequality and lemma~\ref{fo} is proved.\hfill$\square$

 \subsection{ Quantization of pseudo-holomorphic curves on  degenerating Riemann surfaces}
 
We consider a closed Riemann surfaces $(\Sigma, h)$, where $\Sigma$ is smooth compact surface without boundary, and $h$ a metric on $\Sigma$. Since we are only interested in the conformal structure of $\Sigma$, we can assume, thanks to the uniformization theorem, see \cite{Hub},  that $h$ has constant scalar curvature. We consider $(N,J)$ to be a smooth almost-complex manifold and we look at pseudo-holomorphic curves between $(\Sigma,h)$ and $(N,J)$, in other words we consider applications $u \in W^{1,2}(\Sigma, N)$ satisfying
\beq
\label{y7a}
\frac{\p u}{\p x} = J(u)  \frac{\p u}{\p y} \quad,
\eeq
where $z=x+iy$ are some local conformal coordinates on $\Sigma$. These objects are fundamental in symplectic geometry, see \cite{McDS}. In the study of the {\it moduli space}
of pseudo-holomorphic curves in an almost complex manifold, the compactification question comes naturally. In other words it is of first importance to understand and describe
how sequences of pseudo-holomorphic curves with possibly degenerating conformal class behave at the limit. \\

The so-called Gromov's compactness theorem \cite{Gr}, see also \cite{PW}, \cite{Si} and \cite{Hum}, provides an answer to this question.

\begin{thm}
\label{gr}\cite{Gr}
Let $(N,J)$ a compact almost manifold, $\Sigma$ a closed surface and $(j_n)$ a sequence of complex structures on $\Sigma$. Assume $u_n : (\Sigma,j_n)\rightarrow (N,J)$ is a sequence of pseudo-holomorphic curves of bounded area with respect to an arbitrary metric on $N$. Then $u_n$ converge weakly to some cusp curve\footnote{we refer to chapter 5 of  \cite{Hum} for precise definitions} $\overline{u}: \overline{\Sigma}\rightarrow (N,J)$ and there exists finitely many bubbles, holomorphic maps $(\om^i)_{i=1\cdots l}$ from ${S}^2$ into $(N,J)$, such that, modulo extraction
of a subsequence
\[
\lim_{n\rightarrow+\infty} E(u_n)=E(\overline{u})+\sum_{i=1}^lE(\omega^i)\quad.
\]
\hfill $\Box$
\end{thm}

In fact the bound, on the energy is not necessary assuming that the target manifold is symplectic, i.e if there is $\omega$ a closed $2$-form on $N$ compatible with $J$. In deed, in that case, see chapter 2 of \cite{McDS} for instance,  all $u: \Sigma \rightarrow N(J,\omega)$, regular enough, satisfies
$$A(u)=\int_\Sigma dvol_{u^* g} \geq \int_\Sigma u^* \omega$$
where $g=\omega(.,J.)$, with equality if and only if $u$ is pseudo-holomorphic. Hence, for symplectic manifold, pseudo-holomorphic curves are area minimizing in their homology class. In particular, they are minimal surfaces, i.e.  conformal and harmonic, and we can use the general theory of harmonic maps, see  remark 4.2 of \cite{Zhu}.\\ 

We propose below a proof of theorem~\ref{gr} that follows the main lines of the most classical one (i.e.  we shall decompose our curves in thin and thick parts at the limit) but 
the argument we provide in order to prove that there is no energy  in the neck and collar regions is new.  We don't make use of the standard isoperimetric machinery but we simply apply the first order Wente's estimate on annuli given by  lemma~\ref{fo} which fits in an optimal way the particular structure of the pseudo-holomorphic equation (\ref{y7a}) .\\
 
\noindent{\it Proof of theorem~\ref{gr}:}\\

The proof consists in splitting the surface in several pieces where the sequence converges either strongly to a non-constant limiting map or weakly to a constant. Then in a second step, we  prove that there is in fact  no energy  in the pieces where the converge is weak.
Note that in contrast to the previous section, in the present case the complex structure of the surface is not fixed and is {\it a priori} free to degenerate.\\ 
Ou aim is to show how lemma~\ref{fo} can be used in this context and therefore we shall be more brief on the classical parts such as the limiting Deligne-Mumford {\it thin-thick} decomposition which is described for instance in \cite{Hum} or in \cite{Zhu}. Observe that. due to the structure of the equation the $\eps$-regularity theorem for pseudo-holomorphic curves is a consequence of theorem~\ref{ereg}    .\\

For simplicity, we will also assume that we have a surface of genus $g$ greater or equal to $2$. Hence let $h_n$ be the hyperbolic metric of volume 1 associated to the complex structure $j_n$,\\

According to the Deligne-Mumford compactification of Riemann surfaces, see chapter 4 of \cite{Hum}, modulo extraction of a subsequence,
$(\Sigma,h_n)$ converges to an hyperbolic Riemannian $(\Sigma,h)$ surface by collapsing $p$ ($0\leq p\leq 3g-3$) pairwise disjoint simple closed geodesics $(\gamma_n^i)$.\\

{\bf Far from the collapsing geodesics}, the metric uniformly converges, and we have a classical "bubble-tree" decomposition, that is to say $u_n$ converges to a pseudo-holomorphic curves  of the $(\Sigma,h)$ expect possibly at finitely many  points where, as in the previous section, $u_n$ is forming  bubbles (i.e  pseudo-holomorphic curves from $\C$ to $N$) 
which are ''connected'' to each other by some {\bf neck regions} $N_n^i=B(a_n^i,\mu_n^i)\setminus B(a_n^i,\lambda_n^i)$ where the weak $L^2$ energy goes to zero, 
\[
\lim_{r\rightarrow 0}\lim_{n\rightarrow +\infty} \Vert \nabla u_n\Vert_{L^{2,\infty}(N_n^i(r))}=0,
\]
where $N_n^i(r)=B(a_n^i,r\mu_n^i)\setminus B(a_n^i,\frac{\lambda_n^i}{r})$. This can be established by combining the fact that, on such annular regions,
the maximal $L^2$ energy of $\nabla u_n$ on dyadic annuli has to vanish (otherwise we would have another bubble) and the fact that lemma~\ref{L2infini} 
applies to this situation.\\

{\bf Near the collapsing geodesics}, our surface becomes asymptotically isometric to an hyperbolic cylinder of the form

$$
A_l=\left\{ z=re^{i\phi}\in \h : 1\leq r \leq e^l, arctan\left(sinh\right(\frac{l}{2}\left)\right)<\phi<\pi - arctan\left(sinh\right(\frac{l}{2}\left)\right)\right\},
$$
where the geodesic correspond to $\left\{ r e^{i\frac{\pi}{2}}\in \h : 1 \leq r \leq e^l\right\}$ and the line $\{r=1\}$ and $\{ r=e^l\}$ are identified via $z\mapsto e^l z$. This is the {\bf collar region}. It is sometimes easier to consider the following  cylindrical parametrization, i.e.

$$P_l=\left\{ (t,\theta) : \frac{2\pi}{l}arctan\left(sinh\right(\frac{l}{2}\left)\right)<t<\frac{2\pi}{l} \left(\pi - arctan\left(sinh\right(\frac{l}{2}\left)\right)\right), 0\leq \theta \leq 2\pi\right\}$$
in this parametrization the constant scalar curvature metric reads
$$
ds^2=\left(\frac{l}{2\pi sin(\frac{lt}{2\pi}) }\right)^2 (dt^2 +d\theta^2),
$$
where the geodesic corresponds to $\{ t=\frac{\pi^2}{l}\}$ and the line $\{\theta=0\}$ and $\{ \theta=2\pi\}$ are identified.

Then, as $l_n$, the length of the degenerating geodesic, goes to zero, $P_{l_n}=[0,T_n]\times S^1$ up to translation, which can be decompose as follows, see proposition 3.1 of \cite{Zhu}.\\

For each such a thin part, one can extracts a subsequence such that the following decomposition holds. There $p\in {\N}$ and exists 2p sequences $(a_n^1)$, $(b_n^1)$, $(a_n^2)$, $(b_n^2)$,\dots, $(a_n^p)$, $(b_n^p)$ of positive numbers between $0$ and $T_n$ such that
$$ \lim_{n\rightarrow +\infty} \frac{b^i_n-a^i_n}{T_n}=0.$$
and  up to rescaling and identifying $]-\infty,+\infty[\times S^1$ with $\C\setminus \{0\}$, there exists a bubble $\omega^i$ (i.e  pseudo-holomorphic curve from $\C$ to $N$) such that

$$u^n\left(\frac{a_n^i+b_n^i}{2} + \frac{t}{b_n^i-a_n^i},\theta \right) \rightarrow \omega^i \hbox{  on  } C^2_{loc}(\C\setminus \{0\})\quad.$$
Moreover, for any $\eps>0$, there exists $r>0$ such that for any $T\in [b_n^i+r^{-1},a_n^{i+1}-r^{-1}]$
\beq
\label{aq1}
\int_{[T,T+1]\times S^1}|\nabla u_n|^2\le \eps\quad.
\eeq
Denoting $J_n^i=[a_n^i,b_n^i]\times S^1$, $I_n^0=[0,a_n^1]\times S^1$, $I_n^i=[b_n^i,a_n^{i+1}]\times S^1$ and $I_n^p=[b_n^p,T_n]\times S^1$
and $I_n^i(r)=[b_n^i+r^{-1},a_n^{i+1}-r^{-1}]$, (\ref{aq1}) combined with lemma~\ref{L2infini}  implies  
\beq
\label{aq2}
\lim_{r\rightarrow 0}\lim_{n\rightarrow +\infty} \Vert \nabla u_n\Vert_{L^{2,\infty}(I_n^i(r))}=0\quad.
\eeq
This decomposition is illustrated by the following picture.

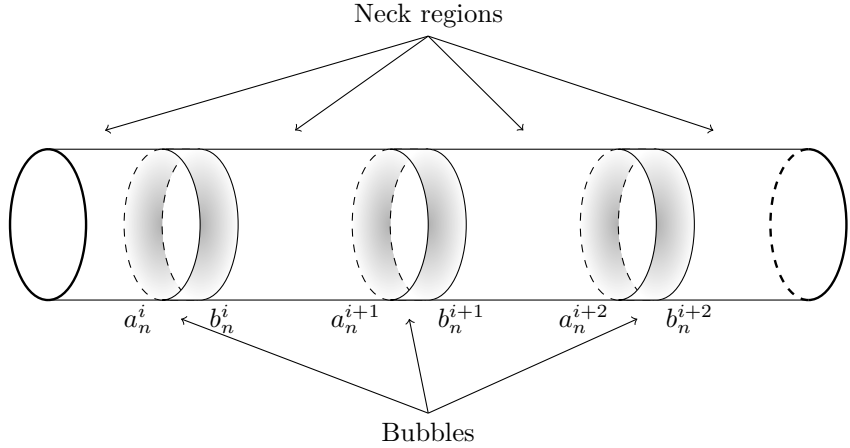
\begin{figure}[!h]
\centering
\input{laurain-fig2.tex}
\caption{Decomposition in necks and bubbles}
\end{figure}

As in the previous section, in order to prove that there is no energy at the limit in the {\bf neck regions} of the thin parts, we combine
the vanishing of the $L^{2,\infty}-$norm given by (\ref{aq2}) with a uniform estimate on the $L^{2,1}$ norm of $\vert \nabla u^n\vert$ on each $I_n^i(r)$, which is
 a direct consequence of the lemma \ref{fo} applied to the pseudo-holomorphic equation
 $$
 \nabla u_n =J(u_n) \nabla^\bot u_n\quad.
 $$ 
 This concludes the proof of theorem~\ref{gr}.\hfill$\square$ 

\begin{rem}
Here again, in addition to the fact that  our argument is not specific to $J$-holomorphic curves, our proof, in comparison with previous ones
such as the one given in \cite{Zhu}, has the advantage to require less regularity on the target manifold $N$. In fact, following the approach of \cite{Pa} or \cite{LiWa}, in order to establish the angular energy quantization, M.Zhu goes through a lower estimate of the following second derivative
$$\frac{d^2}{d\theta^2}\int_{S^1\times\{t\}} \vert u_\theta\vert^2 \, d\theta$$
Such an estimate requires for the metric of $N$ to be at least $C^2$. In the alternative proof we are providing, in order to apply lemma \ref{fo}, we only require the almost complex structure and the compatible metric to be $C^1$ which corresponds to a weakening of the assumption of magnitude 1 in the derivative.
\end{rem}

\subsection{ Quantification for harmonic maps on a degenerating surfaces, a cohomoligical condition.}
The aim of this section is to shed a new light on the quantization for  harmonic maps on a degenerating surfaces, which has been fully described by M.Zhu in \cite{Zhu}.

The main result in the present subsection is the following results which connects energy quantization for harmonic maps into spheres with a cohomological condition.

\begin{thm}
\label{cohom}
Let $(\Sigma,h_n)$ be a sequence of closed Riemann surfaces equipped with their constant scalar curvature metric with volume 1. Let $u_n$ be a sequence
of harmonic maps from $(\Sigma,h_n)$ into the unit sphere $S^{m-1}$ of the euclidian space ${\R}^m$. Assume
\[
\limsup_{n\rightarrow +\infty} E(u_n)<+\infty
\]
and assume that the following closed forms
\[
\forall i,j=1\cdots m\quad\quad\star(u_n^i\ du_n^j-u_n^j\ du_n^i)
\]
are all exact. Then the energy quantization holds : modulo extraction  of  a subsequence, on each component of the limiting thick part,
$u_n$ converges strongly, away from the punctures, to some limiting harmonic map $u$ and  there exists finitely many bubbles, holomorphic maps $(\om^i)_{i=1\cdots l}$ from ${S}^2$ into $(N,J)$, - forming possibly  both on the thick and the thin parts - such that, modulo extraction
of a subsequence
\beq
\label{co1}
\lim_{n\rightarrow+\infty} E(u_n)=E(u)+\sum_{i=1}^lE(\omega^i)\quad.
\eeq
.\hfill$\Box$
\end{thm}

\noindent{\it Proof of theorem~\ref{cohom}:}\\

In fact, assuming that our sequence of harmonic maps $u_n$ get valued into a sphere $S^m$, the equation simply write
$$ \Delta u^i_n = \left( u_n^i \nabla (u_n)_j-(u_n)_j \nabla u^i_n\right)\nabla u_n^j.$$
But $div\left( u_n^i \nabla (u_n)_j-(u_n)_j \nabla u^i_n\right)=0=d(*u_n\wedge du_n)$. Hence {\bf assuming that the closed $\wedge^2{\R}^m$ valued 1-form $\star(u_n\wedge du_n)$ is exact }, there exists $b_n\in W^{1,2}$ such that 
$$\star(u_n\wedge du_n)=db_n ,$$
and 
$$\Vert b_n \Vert_{W^{1,2}} =O\left( \Vert u_n \Vert_{W^{1,2}}\right).$$
Then we have
$$div(\nabla u_n - \nabla^\bot b_n\ u_n)=0.$$
If we are on a neck region such as $B_1\setminus D(0,\eps_n)$, it can be integrated as 
\beq
\label{eqhs}
\nabla u_n = \nabla^\bot b_n\ u_n + \nabla^\bot  c_n + d_n \nabla Log(\rho) ,
\eeq
where $c_n\in W^{1,2}(B_1)$ and $d_n \in\R$. Then we try to control the gradient of the logarithmic part, remarking that 
\be
\begin{split}
\frac{d}{d\rho} \int_0^{2\pi} u_n \, d\theta &= \int_0^{2\pi} \frac{1}{\rho} \frac{\p b_n}{\p \theta}\ u_k \, d\theta + 2\pi\, \frac{d_n}{\rho}\\
&= \int_0^{2\pi} \frac{1}{\rho} \frac{\p b_n}{\p \theta}\ (u_k -\overline{u}^\rho_n) \, d\theta + 2\pi\,  \frac{d_n}{\rho},
\end{split}
\ee
where $\overline{u}^\rho_n$ is the mean value of $u_n$ over $\p B_{\rho}$. Integrating the previous identity from $\eps_n$ to an arbitrary $\rho$ gives
\beq
\label{hs1}
 2\pi (\overline{u}^\rho_n - \overline{u}^{\eps_n}_n)= \int_{\eps_n}^\rho \int_0^\pi \frac{1}{t} \frac{\p b_n}{\p \theta}\ (u_k -\overline{u}^t_n)\, d\theta \, dt + 2\pi\, Log\left(\frac{\rho}{\eps_n}\right).
\eeq
And, thanks to Poincar\'e's inequality, we get 
\beq
\label{hs2}
\left\vert  \int_{\eps_n}^\rho \int_0^\pi \frac{1}{t } \frac{\p b_n}{\p \theta}\ (u_k -\overline{u}^t_n)\, d\theta \, dt\right\vert \leq \Vert \nabla b_n\Vert_2  \Vert \nabla u_n\Vert_2.
 \eeq
Then, combining (\ref{hs1}) and (\ref{hs2}), we finally obtain that 
$$d_n = O \left(\frac{1}{Log\left(\frac{1}{\eps_n}\right)}\right) .$$
Which implies, as in the proof of  lemma~\ref{fo},  that the $L^{2,1}$-norm of $d_n \nabla Log(\rho)$ in $B_1\setminus D(0,\eps_n)$ is uniformly bounded. By the mean of  lemma \ref{fo} and thanks to (\ref{eqhs}), we see that  he $L^{2,1}$-norm of $\nabla (u_n -d_n Log(\rho))$ is also uniformly bounded and these two uniform bounds imply the uniform $L^{2,1}$ bound
of $\nabla u_n$ in neck regions. Combining the  uniform $L^{2,1}$ bound
of $\nabla u_n$ in neck regions together with the lemma~\ref{L2infini} gives the desired energy quantization (\ref{co1}) and theorem~\ref{cohom} is proved.\hfill $\Box$ \\

More generally we can raise the following question :
{\bf Considering  a sequence of harmonic maps from a degenerating surface to a general target manifolds, is there is  a simple cohomological condition similar as the one in theorem~\ref{cohom} ensuring the quantization of the energy in collar region?}

\subsection{Energy Quantization for Willmore Surfaces.}

Finally we would like to recall a last application of lemma~\ref{fo} that has been used in a recent work by Y.Bernard and T.Riviere in \cite{BR}
for proving Energy Quantization for sequences of Willmore surfaces with uniformly bounded energy and non-degenrating conformal classes. The problem can be described as follows : 
for a sufficiently smooth immersion $u : \Sigma \rightarrow \R^m$, where $\Sigma$ is a closed two dimensional Riemannian surface, we can define its mean curvature vector $\vec{H}$ and we consider  the following functional
$$W(u)= \int_\Sigma  |\vec{H}|^2\, u^*(dy).$$
where $u^*(dy)$ denotes the metric induced on $\Sigma$ by the immersion $u$. This functional is called, the Willmore functional and is known to be  conformally invariant (see \cite{Ri6}). Critical points to the functional $W$ are called {\it Willmore immersions} or 
{\it Willmore surfaces}. Hence as for harmonic maps or pseudo holomorphic curves the question of the quantization of sequences of Willmore surfaces arise naturally. The second author has developed appropriate tools to study weak critical points to $W$ in \cite{Ri4} and \cite{Ri5} and proved the epsilon-regularity for these weak critical points.
Using in particular lemma~\ref{fo} the following energy quantization has been established
\begin{thm} \cite{BR} 
Let $u_n$ be a sequence of Willmore immersions of a closed surface $\Sigma$. Assume that
$$\limsup_{n \rightarrow +\infty } W(u_n) <+\infty $$
that the conformal class of $u_n^*(\xi_{\R^m})$ remains within a compact subdomain of the moduli space of $\Sigma$. Then,
modulo extraction of a subsequence, the following energy identity holds
$$\lim_{n \rightarrow +\infty } W(u_n) = W(u_\infty)+\sum_{l=1}^L W(\omega_l)  + \sum_{k=1}^K (W(\Omega_k) - 4\pi \theta_k)  $$
where $u_\infty$ is a possibly branched smooth Willmore immersion of $\Sigma$. The maps $\omega_l$ and
$\Omega_k$ are smooth , possibly branched, Willmore immersions of $S^2$ and $\theta_k$ is the integer density of the current $(\Omega_k)_*(S^2)$ at some point $p_k \in \Omega_k(S^2)$, namely
$$\theta_k= \lim_{\rho \rightarrow 0}\frac{\mathcal{H}^2(B_r(p)\cap \Omega_k(S^2))}{\pi \rho^2}.$$
\end{thm}
%%%%%%%%%%%%%%%%%%%%%%%%%%%%%%%
\appendix
\section{Lorentz Estimates on Harmonic Functions.}
Here we prove two lemmas on harmonic which insure that we can control the $L^{2,1}$-norm by the $L^{2}$ on a smaller domain up to some appropriate boundary condition.

\begin{lemma}
\label{l1}
Let $0<\eps<\frac{1}{2}$ and $f:B_1\setminus B_\eps \rightarrow \R$ an harmonic function which satisfies
\beq
\label{bord}
\begin{split}
&f=0 \hbox{ on } \partial B_1,\\
&\int_{\partial B_\eps} f\;d\sigma =0.
\end{split}
\eeq
Then, for all $\lambda>1$, there exists positive a constant  $C(\lambda)$ independent of $\eps$ and $f$ such that 
$$ \Vert \nabla f\Vert_{L^{2,1}(B_1\setminus B_{\lambda\eps})} \leq C(\lambda) \Vert \nabla f\Vert_2 .$$
\hfill$\Box$
\end{lemma}
\noindent{\it Proof of lemma~\ref{l1} :}\\

We start by decomposing $f$ as a Fourier series, which gives
$$f(\rho,\theta) = c_0 + d_0\ln(\rho) +\sum_{n\in \Z^{*}} (c_n\rho^n+d_n \rho^{-n}) e^{in\theta}.$$
Hence, using (\ref{bord}), we easily proved that $c_0=d_0=c_n+d_n=0$, then we get 
$$f(\rho,\theta) = \sum_{n\in \Z^{*}} c_n(\rho^n- \rho^{-n}) e^{in\theta}.$$
Then we estimate the gradient as follows
\be
\vert \nabla f(\rho,\theta) \vert \leq 2 \sum_{n\in \Z^{*}} \vert n\, c_n\vert (\rho^{n-1} + \rho^{-n-1}).
\ee
Then, we estimate  the $L^{2,1}$-norm of the $f_m(z)= \vert z\vert^m$ on $B_1\setminus B_{\lambda\eps}$, for $m\in \Z\setminus \{-1\}$ and $\lambda\in]1,2]$, which gives
\beq
\begin{split}
&\Vert f_m\Vert_{L^{2,1}(B_1\setminus B_{\lambda\eps})} \leq \sqrt{\pi} \int_{0}^{(\lambda\eps)^m} t^{\frac{1}{m}}\; dt \leq 2\sqrt{\pi} (\lambda\eps)^{m+1} \hbox{ for } m<-1\\
&\hbox{ and } \\
&\Vert f_m\Vert_{L^{2,1}(B_1\setminus B_{\lambda\eps})}\leq \sqrt{\pi}\hbox{ for } m\geq 0.
\end{split}
\eeq
Here we use the following characterization (\ref{d1}).
Hence we get 
\be
\Vert \nabla f \Vert_{L^{2,1}(B_1\setminus B_{\lambda\eps})} \leq 4\sqrt{\pi}\left( \sum_{n>0} \vert n\, c_n\vert \left( (\lambda\eps)^{-n} +1\right) + \sum_{n<0} \vert n\, c_n\vert \left( (\lambda\eps)^{n} +1\right) \right).
\ee
Hence, thanks to the Cauchy-Scharwz and the fact that $\lambda>1$, we get 
\be
\Vert \nabla f \Vert_{L^{2,1}(B_1\setminus B_{\lambda\eps})} \leq 8\sqrt{\pi} \left(\sum_{n\not= 0}  \vert n\vert \lambda^{-2\vert n\vert} \right)\left(\sum_{n\not= 0} \vert n\vert \,\vert c_n\vert^2 \eps^{-2\vert n\vert}\right)^\frac{1}{2} .
\ee
Finally we compute the $L^2$-norm of $\nabla f$
$$\Vert \nabla f\Vert_2 =\left(2\pi \int_{\eps}^{1} \sum_{n\not=0} \vert n\, c_n\vert^2 (\rho^{2n-2} + \rho^{-2n-2})\,\rho d\rho\right)^\frac{1}{2} \geq  \sqrt{\frac{\pi}{2}}\left(\sum_{n\not=0 } \vert n\vert \, \vert c_n\vert^2 \eps^{-2\vert n\vert}\right)^\frac{1}{2} $$
which achieves the proof of  lemma~\ref{l1} .\hfill$\square$\\

\begin{lemma}
\label{l3}
 let $0<\eps<\frac{1}{4}$ and $f:B_1\setminus B_\eps \rightarrow \R$ an harmonic function which satisfies
\beq
\label{bord2}
\begin{split}
&\int_{\partial B_\eps} f\;d\sigma =0,\\
&\left\vert \int_{\partial B_1 } f  \, d\sigma \right\vert \leq K,\\
\end{split}
\eeq
where $K$ is a constant independent of $\eps$. Then, for all $0 <\lambda<1$ there exists positive constant  $C(\lambda)$ independent of $\eps$ and $f$ such that 

$$ \Vert \nabla f\Vert_{L^{2,1}( B_{\lambda}\setminus   B_{\lambda^{-1}\eps })} \leq C(\lambda) (\Vert \nabla f\Vert_2+1) .$$
\hfill $\Box$
\end{lemma}

{\it Proof of lemma~\ref{l3} :}\\

We start by decomposing $f$ as a Fourier series, which gives
$$f(\rho,\theta) = c_0 + d_0\ln(\rho) +\sum_{n\in \Z^{*}} (c_n\rho^n+d_n \rho^{-n}) e^{in\theta}.$$
Hence, using (\ref{bord2}), we easily proved that $c_0+d_0\,ln(\eps)=0$ and $\vert c_0 \vert=O(1)$. Hence 
\beq
\label{est}
d_0=O\left( \frac{-1}{ln(\eps)}\right).
\eeq
Then we estimate the gradient as follows
\be
\vert \nabla f(\rho,\theta) \vert \leq \vert d_0\vert \frac{1}{\rho} + \sum_{n\in \Z^{*}} \vert n\, c_n\vert \rho^{n-1} +\vert n\, d_n\vert  \rho^{-n-1}.
\ee
Then, we estimate  the $L^{2,1}$-norm of the $f_m(z)= \vert z\vert^m$ on $ B_{\lambda}\setminus   B_{\lambda^{-1}\eps }$, for $m\in \Z\setminus \{-1\}$ and $0<\lambda<1$, which gives
\beq
\label{est2}
\begin{split}
&\Vert f_m\Vert_{2,1}\leq \sqrt{\pi} \int_{0}^{(\lambda^{-1}\eps)^m} t^{\frac{1}{m}}\; dt \leq 2\sqrt{\pi} (\lambda^{-1}\eps)^{m+1} \hbox{ for } m<-1\\
&\Vert f_m\Vert_{2,1}\leq \sqrt{\pi }\lambda^m\hbox{ for } m\geq 0,\\
&\hbox{ and } \\
& \Vert f_{-1}\Vert_{2,1}=O(-\log(\eps)).
\end{split}
\eeq
Here we use the following characterization (\ref{d1}). Thanks to (\ref{est}) and (\ref{est2}), we get 
\be
\begin{split}
\Vert \nabla f \Vert_{L^{2,1}(B_{\lambda}\setminus   B_{\lambda^{-1}\eps })} &\leq 2\sqrt{\pi}\left( \sum_{n>0}  \left(\vert n\, c_n\vert \lambda^{n} +\vert n\, d_n\vert  (\lambda^{-1}\eps)^{-n}\right) + \sum_{n<0}\left(\vert n\, c_n\vert (\lambda^{-1}\eps)^{n} +\vert n\, d_n\vert \lambda^{-n}\right) \right) \\
&+O(1).
\end{split}
\ee
Hence, thanks to the Cauchy-Scharwz and the fact that $0<\lambda<1$, we get 
\be
\begin{split}
\Vert \nabla f \Vert_{L^{2,1}(B_{\lambda}\setminus   B_{\lambda^{-1}\eps })} &\leq 4\sqrt{\pi} \left(\sum_{n\not= 0}  \vert n\vert \lambda^{2\vert n\vert} \right) \left(\sum_{n< 0} \vert n\vert \,(\vert c_n\vert^2+\vert d_{-n}\vert^2) \eps^{-2\vert n\vert}+
\sum_{n>0} \vert n\vert \,(\vert c_n\vert^2+\vert d_{-n}\vert^2) 2^{- n}\right)^\frac{1}{2}\\
& +O(1).
\end{split}
\ee
Finally we compute the $L^2$-norm of $\nabla f$
\be
\begin{split}
\Vert \nabla f\Vert_2 &= \vert d_0\vert  \left( \int_\eps^1 \frac{1}{\rho}\, d\rho\right)^\frac{1}{2} + \left(2\pi \int_{\eps}^{1} \sum_{n\not=0} \left(\vert n\, c_n\vert^2 \rho^{2n-2} +\vert n\, d_n\vert^2  \rho^{-2n-2}\right)\,\rho d\rho\right)^\frac{1}{2}\\
& \geq  \sqrt{\frac{\pi}{2}}\left(\sum_{n<0 } \vert n\vert \, (\vert c_n\vert^2 + \vert d_{-n}\vert^2) \eps^{-2\vert n\vert}+
\sum_{n>0} \vert n\vert \,(\vert c_n\vert^2+\vert d_{-n}\vert^2) 2^{- n}\right)^\frac{1}{2} 
\end{split}
\ee
which achieves the proof of lemma~\ref{l3}.\hfill$\square$\\

\end{document}

%% file: laurain-fig1.tex
\begin{tikzpicture}[scale=0.25]
\draw (-12,0).. controls (-12,-2) and (-10,-5).. (-6,-5) .. controls (-4,-5) and (-2,-2) ..  (0,-2).. controls (2,-2) and (4,-3) .. (6,-3) .. controls (7,-3) and (10,-2).. (10,0)..controls (10,2) and (6,4) .. (4,3)..controls (2,2) and (0,1) .. (-2,2).. controls (-4,3) and (-6,4) .. (-8,4).. controls (-10,4) and (-12,2)..(-12,0);

\draw (-12,0).. controls (-12,-2) and (-10,-5).. (-6,-5) .. controls (-4,-5) and (-2,-2) ..  (0,-2).. controls (2,-2) and (4,-3) .. (6,-3) .. controls (7,-3) and (10,-2).. (10,0)..controls (10,2) and (6,4) .. (4,3)..controls (2,2) and (0,1) .. (-2,2).. controls (-4,3) and (-6,4) .. (-8,4).. controls (-10,4) and (-12,2)..(-12,0);

\draw (-9,0) .. controls (-8,-2) and (-6,-2)..(-5,0);
\draw(-8.7,-0.5) .. controls (-8,1) and (-7,1)..(-5.3,-0.5);

\draw (5,0) .. controls (6,-1) and (7,-1)..(8,0);
\draw (5.3,-0.3) .. controls (6,1) and (7,1)..(7.7,-0.3);

\draw (-4.7,-2.4) node {$\times$};
\draw (0.7,-0.4) node {$\times$};
\draw (4,-1.4) node {$\times$};
\draw (-5.7,2.4) node {$\times$};

\draw (-4.7,-2.4)  circle (20pt) ;

\draw[->] (-4.7,-3.1) .. controls  (-4.7,-4.4) and (-4.7,-5.4) .. (-4.7,-6.4);

\draw (-5.7,-7.5) node {$\times$};
\draw (-3.5,-9.5) node {$\times$};
\draw (-5,-9) node {$\times$};

\draw (-4.7,-9)  circle (74pt) ;
\draw (-5,-9)  circle (20pt) ;

\draw[->, dashed] (-5,-9.7) .. controls  (-5,-11) and (-5,-12) .. (-5,-13);
\end{tikzpicture}

%% file: laurain-fig2.tex
\begin{tikzpicture}[scale=0.5]
\draw (-10,2)--(10,2); 
\draw (-10,-2)--(10,-2);
\draw[line width=0.9pt] (-10,2) arc (90:270: 1cm and 2cm);
\draw[line width=0.9pt] (-10,2) arc (90:-90: 1cm and 2cm);

\draw[line width=0.9pt] (10,2) arc (90:-90: 1cm and 2cm);
\draw[dashed, line width=0.9pt] (10,2) arc (90:270: 1cm and 2cm);

\draw [draw=black,dashed, inner color=gray!60]  (-7,2)
arc (90:270:1cm and 2cm) --(-6,-2)
arc (-90:-270:1cm and 2cm ) -- cycle;

\draw [draw=black, inner color=gray!60]  (-7,2) -- (-6,2)
arc (90:-90:1cm and 2cm) -- (-7,-2)
arc (-90:90:1cm and 2cm ) -- cycle;

\draw [draw=black,dashed, inner color=gray!60]  (-1,2)
arc (90:270:1cm and 2cm) --(-0,-2)
arc (-90:-270:1cm and 2cm ) -- cycle;

\draw [draw=black, inner color=gray!60]  (-1,2) -- (-0,2)
arc (90:-90:1cm and 2cm) -- (-1,-2)
arc (-90:90:1cm and 2cm ) -- cycle;

\draw [draw=black,dashed, inner color=gray!60]  (5,2)
arc (90:270:1cm and 2cm) --(6,-2)
arc (-90:-270:1cm and 2cm ) -- cycle;

\draw [draw=black, inner color=gray!60]  (5,2) -- (6,2)
arc (90:-90:1cm and 2cm) -- (5,-2)
arc (-90:90:1cm and 2cm ) -- cycle;

\draw (-7,-2.5) node[left] {$a_n^i$}; 
\draw (-6,-2.5) node[right] {$b_n^i$}; 

\draw (-1,-2.5) node[left] {$a_n^{i+1}$}; 
\draw (0,-2.5) node[right] {$b_n^{i+1}$}; 

\draw (5,-2.5) node[left] {$a_n^{i+2}$}; 
\draw (6,-2.5) node[right] {$b_n^{i+2}$}; 

\draw (0,-5) node[below] {Bubbles};
\draw[->] (0,-5) -- (-6.5,-2.5);
\draw[->] (0,-5) -- (-0.5,-2.5);
\draw[->] (0,-5) -- (5.5,-2.5);

\draw (0,5) node[above] {Neck regions};
\draw[->] (0,5) -- (-8.5,2.5);
\draw[->] (0,5) -- (-3.5,2.5);
\draw[->] (0,5) -- (2.5,2.5);
\draw[->] (0,5) -- (7.5,2.5);

\end{tikzpicture}